\documentclass{article}
\pdfoutput=1
\usepackage{graphicx, amssymb, amsmath}

\usepackage{amsthm}
\usepackage{mathrsfs}
\usepackage{color}
\usepackage{graphicx}
\usepackage{caption}
\usepackage[ruled]{algorithm2e}  
\usepackage{amsmath}

\usepackage{upgreek}
\usepackage{subcaption}
\usepackage{color,    soul}
\usepackage{cite} 
\usepackage[colorlinks=true, linkcolor=black, citecolor=blue, urlcolor=blue]{hyperref}
\usepackage{setspace}
\sethlcolor{blue}
\usepackage{float}
\usepackage{tabularx}
\usepackage{multirow}
\usepackage{makecell}
\usepackage{booktabs}

\usepackage{amsmath}
\usepackage{multirow}
\usepackage{diagbox}
\usepackage{subcaption}
\usepackage{chngcntr}
\counterwithin{figure}{section}  
\counterwithin{table}{section}   
\usepackage{subcaption}

\def\k{\kern .5em}
\def\er{\kern .2em}

\usepackage{tikz}
\usetikzlibrary{positioning, shapes.geometric}

\begin{document}
	
	\counterwithin{algocf}{section} 
	
	\date{}
	\author{}
	\newcommand{\be}{\begin{equation}}
		\newcommand{\ee}{\end{equation}}
	\newcommand{\ba}{\begin{array}}
		\newcommand{\ea}{\end{array}}
	
	\newcommand{\beas}{\begin{eqnarray*}}
		\newcommand{\eeas}{\end{eqnarray*}}
	\newcommand{\bea}{\begin{eqnarray}}
		\newcommand{\eea}{\end{eqnarray}}
	\newcommand{\ome}{\Omega}

	\newtheorem{theorem}{Theorem}[section]
	\newtheorem{lemma}{Lemma}[section]
	\newtheorem{remark}{Remark}[section]
	\newtheorem{proposition}{Proposition}[section]
	\newtheorem{definition}{Definition}[section]
	\newtheorem{corollary}{Corollary}[section]

\newtheorem{theo}{Theorem}[section]
\newtheorem{lemm}{Lemma}[section]
\newcommand{\blem}{\begin{lemma}}
	\newcommand{\elem}{\end{lemma}}
\newcommand{\bthe}{\begin{theorem}}
	\newcommand{\ethe}{\end{theorem}}
\newtheorem{prop}{Proposition}[section]
\newcommand{\bprop}{\begin{proposition}}
	\newcommand{\eprop}{\end{proposition}}
\newtheorem{defi}{Definition}[section]
\newtheorem{coro}{Corollary}[section]
\newtheorem{algo}{Algorithm}[section]
\newtheorem{rema}{Remark}[section]
\newtheorem{property}{Property}[section]
\newtheorem{assu}{Assumption}[section]
\newtheorem{exam}{Example}[section]

\renewcommand{\theequation}{\arabic{section}.\arabic{equation}}
\renewcommand{\thetheorem}{\arabic{section}.\arabic{theorem}}
\renewcommand{\thelemma}{\arabic{section}.\arabic{lemma}}
\renewcommand{\theproposition}{\arabic{section}.\arabic{proposition}}
\renewcommand{\thedefinition}{\arabic{section}.\arabic{definition}}
\renewcommand{\thecorollary}{\arabic{section}.\arabic{corollary}}
\newcommand{\lan}{\langle}
\newcommand{\curl}{{\bf curl \;}}
\newcommand{\rot}{{\rm curl}}
\newcommand{\grad}{{\bf grad \;}}
\newcommand{\dvg}{{\rm div \,    }}
\newcommand{\ran}{\rangle}
\newcommand{\bR}{\mbox{\bf R}}
\newcommand{\bRn}{{\bf R}^3}
\newcommand{\Coinf}{C_0^{\infty}}
\newcommand{\disp}{\displaystyle}
\newcommand{\ra}{\rightarrow}
\newcommand{\Ra}{\Rightarrow}
\newcommand{\ud}{u_{\delta}}
\newcommand{\Ed}{E_{\delta}}
\newcommand{\Hd}{H_{\delta}}
\newcommand\varep{\varepsilon}
\newcommand{\RNum}[1]{\uppercase\expandafter{\romannumeral #1\relax}}
\newcommand{\tabincell}[2]{\begin{tabular}{@{}#1@{}}#2\end{tabular}}
\title{A Novel Iterative Algorithm Incorporating Machine Learning for the Poisson-Nernst-Planck Equations in Three-Dimensional Ion Channel Simulation}

\author{Yuyan Ao$^{1, 2}$ \and Benzhuo Lu $^{3}$\and Shi Shu $^{1}$ \and Ying Yang $^{2, *}$ }
\footnotetext[1]
{School of Mathematics and Computational Science, Xiangtan University,
	Xiangtan, 411105, China}
    

\footnotetext[2]
{ School of Mathematics and Computing Science, Guangxi Colleges and Universities Key Laboratory of Data Analysis and Computation, Center for Applied Mathematics of Guangxi (GUET), Guilin University of Electronic Technology, Guilin, 541004, Guangxi, China. \textbf{*Corresponding author.} E-mail: yangying@lsec.cc.ac.cn}

\footnotetext[3]
{ State Key Laboratory of Scientific and Engineering Computing, Institute of Computational Mathematics and Scientific Engineering Computing, Academy of Mathematics and Systems Science, Chinese Academy of Sciences, Beijing, 100190, China. \textbf{*Corresponding author.} E-mail: bzlu@lsec.cc.ac.cn}

\maketitle

\noindent
{\bf Abstract}\quad 
The Poisson Nernst Planck (PNP) equation is widely used to simulate ion channel transport, but its strong nonlinearity and coupling characteristics make numerical solutions extremely challenging, and decoupling methods are usually required. The Gummel iteration method is widely adopted as a common decoupling algorithm owing to its simplicity, ease of implementation, and applicability even under complex conditions. But its convergence efficiency highly depends on the selection of the relaxation parameter, unsuitable parameters will lead to slow convergence or even divergence. At present, parameter selection relies heavily on manual experience, which is inefficient and unstable, significantly hindering the solving efficiency for practical problems such as ion channel simulations. To address this, we propose a machine learning-based parameter selection method, which can select near-optimal relaxation parameters for the Gummel iteration of PNP equations for ion channels. Compared to manual empirical approaches, this method substantially enhances both the stability and efficiency of the solution process. Notably, in the field of ion channels, the collection of training data for selecting relaxation parameters via machine learning faces two major challenges. First, the complexity of real ion channel structures often leads to convergence issues in numerical simulations, resulting in a scarcity of effective training samples. Second, extracting geometric features of the channels is particularly difficult. To address the first challenge, we propose an artificial channel generation algorithm to construct a suitable training dataset. For the second challenge, inspired by point cloud algorithms from the 3D printing field, we introduce a geometric feature extraction algorithm that effectively captures key structural features and provides critical input for the model. Finally, based on the machine learning-based parameter selection algorithm, we develop a coupled finite element and Gummel iteration solver, and more stable convergence and higher computational efficiency have been achieved in ion channel simulations. To validate the method, we conduct tests on two real channels with complex geometries that differ significantly from the training samples. Numerical results demonstrate that this method has a significant effect on improving the convergence efficiency and stability of Gummel iteration, confirming its applicability in practical ion channel simulations.
\noindent



\section{Introduction}\label{sec1}
\noindent

Ion channels are a type of special protein micropores that penetrate biological membranes. They can open or close under specific conditions, enabling specific ions to selectively exchange across the cell membrane along the electrochemical gradient, thereby participating in the regulation of numerous cellular physiological activities\cite{hille99}. 
The malfunction of ion channels is associated with numerous diseases. 
Ion channels have been demonstrated as valuable targets for the treatment of epilepsy, chronic pain, and other diseases \cite{clare10}.
The PNP equations are a classical model used to describe the convective diffusion process of charge concentration distribution under the action of an electric field. They are proposed by Nernst \cite{nernst1889elektromotorische} and Planck \cite{planck1890ueber}, and have been widely used in the numerical simulations of semiconductor devices \cite{markowich1985stationary, Jerome1996Analysis}, nanopore systems \cite{Daiguji04,xu18}, and biological ion channels \cite{eisenberg1998ionic, saraniti01, zheng11}. 
The PNP equations for simulating ion channels are characterized by strong coupling, nonlinearity, irregular interfaces, and multiple singular points.
Due to the complexity of the practical problem, analytical solutions are difficult to obtain.
 Many numerical methods such as finite difference method\cite{kurnikova1999FDM,im2002FDM,bolintineanu2009FDM,liu2014FDM}, spectral element method\cite{Hollerbach01}, finite volume method \cite{mathur2009FVM}, boundary element method \cite{zhou2008BEM} and finite element method (FEM)\cite{Tu2013Aparallelfinite, tu2015stabilized,Zhang2020,xie2020effective} have been employed for PNP equations in ion channel simulations.

  Among these numerical methods, FEM is often used to attempt to simulate ion channels with irregular structures due to its good ability to handle complex geometries and flexible boundary conditions, thereby generating approximate high-precision PNP solutions \cite{Tu2013Aparallelfinite,tu2015stabilized,Chaudhry15,Zhang2020,xie2020effective,shen24,xie24}. Tu et al. \cite{Tu2013Aparallelfinite} have developed a finite element simulator for ion transport through three-dimensional ion channel systems, and the simulation results for three real ion channel systems agree well with experimental results. 
    Zhang et al. \cite{Zhang2020} have designed three initialization methods to provide good initial \textcolor{black}{guesses} for the PNP system, and effectively applied them to the simulations of two potassium channels. Xie and Lu \cite{xie2020effective} have presented \textcolor{black}{an} effective finite element solver for a PNP ion channel model with periodic boundary value conditions. 
  Recently,  an inverse averaging FEM was developed in \cite{shen24} for solving the size-modified PNP (SMPNP) equations (which consider ion size effect compared with classic PNP equations), and it is shown \textcolor{black}{to} be more effective and robust than the traditional FEM in simulations of biological systems. Xie \cite{xie24} \textcolor{black}{has} presented an efficient finite element iterative method for solving the SMPNP model for ion channels, and applied \textcolor{black}{it} to simulating a practical ion channel
  and two mixture solutions of four ionic species. 
  
    Due to the strong coupling of PNP equations, it is common to combine decoupling methods with FEM to obtain decoupled solutions. The Gummel iteration is the main decoupling method combined with FEM. This method was proposed by Gummel \cite{gu64} in 1964, initially used to solve the drift-diffusion model \textcolor{black}{for} semiconductor devices. Due to its high efficiency and ease of implementation, it was later promoted for solving the PNP equations \cite{luzh11, chen2013surface, liu2025existence}. This method has also been applied to solving practical PNP ion channel models.  
            Yang et al. \cite{yang2020superconvergent} combined the Gummel iteration with a superconvergence technique to solve the nonlinear PNP ion channel model. Zhang et al. \cite{Zhang2020} proposed three initialization methods based on model reduction techniques and integrated them with the Gummel iteration for numerical simulations of three-dimensional ion channels. 
            Wang et al. \cite{Wang2021} have proposed a novel stabilized FEM named SUPG–IP, which is combined with the Gummel iteration, and demonstrates much better performance in positivity preserving and robustness than the standard FEM and SUPG method in the ion channel simulation.

    There is another type of decoupled method for the PNP equations, called two-grid decoupled method \cite{yang2019two-grid, shen2019decoupling}. This method employs the finite element solution on the coarse grid to decouple the coupled equations on the fine grid. It is comparable in solution accuracy to the Gummel iteration, yet achieves a significantly higher efficiency. 
        However, difficulties are encountered
 when applying this method to practical ion channel problems. Its implementation requires two-level nested coarse and fine grids, whereas meshes generated by professional or commercial software often fail to satisfy such requirements.  Newton’s method is also been used for solving PNP equations (cf. \cite{mathur2009FVM,bousquet2018newton}), but it is commonly used to solve the entire system of PNP
equations, not used to decouple the system as the Gummel iteration and two-grid decoupled method. 

Machine learning, particularly deep learning, has been extensively applied to the numerical solution of partial differential equations in diverse fields such as materials science and fluid mechanics. Recently, very few works have emerged that apply machine learning to the solution of PNP equations. Lee et al. \cite{Lee2024} have presented a deep learning solver for an improved one-dimensional PNP ion channel model. 
Huang et al. \cite{huang25} have proposed an enriched physics-informed neural network, to solve time-dependent PNP equations. 
However, existing studies that apply machine learning or deep learning methods to the PNP equations are not developed for practical ion-channel PNP models.

In this paper, we propose a novel iterative method incorporating machine learning for the PNP equations in three-dimensional ion channel simulations, to improve the computational efficiency of classical iterative methods.
The key difference from existing work is that we specifically design the machine learning training procedure for practical ion channels. Furthermore, we apply the proposed method to simulations of three realistic ion channels to verify its effectiveness.
To the best of our knowledge, no machine learning or deep learning approaches have been specifically designed for the PNP problems in practical ion channel applications. The novel iterative method proposed here is obtained by incorporating machine learning techniques into the Gummel iteration. In what follows, we elaborate on the motivation for employing machine learning methods in the Gummel iteration.

 Although Gummel iteration has been widely applied in solving PNP equations, its convergence efficiency is highly sensitive to the selection of the relaxation parameter (cf. Section 3.1).
Improper selection of this parameter can significantly increase the number of iterations, or even lead to non-convergence, which is particularly prominent in problems such as the ion channel PNP model that require extensive fine-grained calculations using extremely fine grids. Such problems already have high numerical computation costs, and each additional iteration further increases the computational burden. 
Currently, the selection of relaxation parameters is mostly based on experience, lacking specificity and stability. Specifically, optimal relaxation parameters for a specific ion channel PNP equation are often difficult to directly apply to other channels; even within the same channel, parameters that perform well for one physical quantity may also cause divergence in the iteration process for another physical quantity. To address the above issues, this paper proposes a relaxation parameter selection method based on machine learning for optimizing the Gummel iteration process in solving the PNP equation of ion channels. This method can adaptively select relaxation parameters for different PNP ion channel models and achieve adaptive parameter adjustment in the iterative solution of different physical quantities. 

It is noted that there are also other fields using machine learning to select relaxation parameters, such as Silva et al. \cite{Silva2021Machinelearning}, who used Random Forest Regression to predict the relaxation parameters of the Picard iteration to accelerate the solution of multiphase porous media flow. Jiang et al. \cite{GPR} introduced Gaussian Process Regression to predict the parameters of the Alternating Direction Implicit Method to improve the efficiency and accuracy of the solution. Mahmoodzadeh et al. in \cite{Mahmoodzadeh2022} also explored the use of Gaussian Process Regression, Support Vector Regression, Decision Trees, and Long Short-Term Memory Networks to indirectly estimate rock strength parameters in a faster and more economical way. Uyan et al. in \cite{uyan2023mechine} carried out research on using machine learning methods to predict the parameters of underwater sensor networks to improve the performance of underwater sensor networks. These studies demonstrate that the combination of machine learning techniques can effectively improve the efficiency and effectiveness of traditional numerical solution methods. 

In this paper, we use machine learning to select relaxation parameters for Gummel iteration, aiming to accelerate the solution of PNP equations.  However, there are two challenges in the process of collecting training sets when using machine learning to select relaxation parameters. 
First, to improve the solution efficiency of iterative methods, high-resolution numerical solutions are important training samples.
 However, the ion channel models of the PNP equations that traditional numerical methods can solve are extremely limited, typically restricted to relatively small channels (with the number usually being single-digit). This results in a severe scarcity of training samples in the form of high-resolution numerical solutions. Second, the primary difference among different channels resides in their geometric shapes, which also constitutes a critical factor influencing numerical convergence \cite{tai2008molecular}. However, the geometric shapes of most channels are highly irregular (cf. Figures \ref{fig:VDAC} and \ref{fig:KcsA} as examples), and their corresponding physical information is excessively voluminous, making it difficult to be directly used as training data. 
How to accurately capture their geometric features without generating an excessively large number of training samples thus poses a key challenge.

%

To address the problem of insufficient training samples, we first proposed an artificial ion channel generation algorithm, through which the number of samples can be enriched. On this basis, we further constructed two algorithms for extracting the geometric characteristics of channels based on the relevant algorithms in point clouds. It should be noted that although our machine learning algorithm uses artificial channels as training data, it is only applied to predicting the relaxation parameters of the iterative method, acting merely as an adaptive parameter generator. Since the core iteration still relies on the conventional Gummel iteration, the use of artificial channels does not introduce any bias into the numerical solutions. 
Finally, we used this series of algorithms to construct the training set jointly. Based on the constructed training set, we have put forward 
 an improved Gummel iteration based on Random 
Forest Regression (RFR) parameter selection.
This algorithm can automatically select relatively optimal relaxation parameters for channels with different shapes and accelerate the solution of the PNP equations in ion channels. In the numerical tests, we selected three typical ion channels with complex structures, the gramicidin A (GA) channel , the voltage dependent anion channel (VDAC) and a realistic potassium channel: the bacterial channel KscA, as test subjects and conducted validation under multiple voltage and ion concentration conditions. The numerical results indicate that the new iterative method incorporating machine learning is more efficient and stable than the classical Gummel iteratation. 
{\color{black} To our knowledge, although deep-learning-based solvers for one‑dimensional simplified PNP ion channel models have been reported \cite{Lee2024}, machine learning approaches for numerical solutions of PNP equations for real, three-dimensional ion channel structures remain underexplored.}

The organization of this paper is as follows: First, Section 2 introduces the PNP equations, weak formulation and the finite element approximation. 
Then, Section 3 first introduces the algorithm of the FEM combining with Gummel iteration for PNP equations (see Algorithm \ref{Gummel}), and the whole solution procedure for PNP ion channel model. After that, we present three algorithms to construct the dataset for training the RFR (see Algorithms \ref{artificial channel generation}-\ref{for Complex Channel}). Based on this, we propose a novel iteration incorporating RFR. 
This method aims to accelerate the convergence speed of the Gummel iteration. In Section 4, we select three ion channels, GA, VDAC and KcsA, as test cases and compare the number of Gummel iterations and computational efficiency between the parameters selected by random forest regression and traditional fixed parameters under specific concentration and voltage conditions. Finally, the paper concludes with a summary section.


\setcounter{lemma}{0}
\setcounter{theorem}{0}
\setcounter{corollary}{0}
\setcounter{equation}{0}
\setcounter{remark}{0}
\section{PNP equations for ion channels}
\begin{figure}[htbp]
\centering
\includegraphics[height=8cm,    width=5cm]{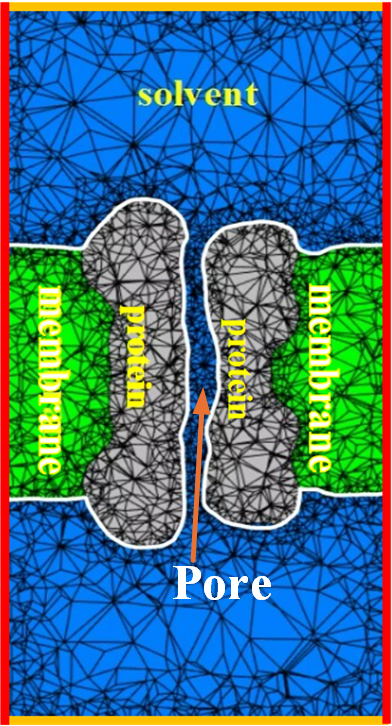}
\caption{2D cut plane $y=0$ for the mesh of GA single channel}
\label{Computationaldomain}
\end{figure}

Let $\Omega \subset R^3$ denote an open domain, $\Omega=\Omega_{m} \cup \Omega_{s}, \Omega_{m} \cap \Omega_{s}=\varnothing$. $\Omega_{m}$ represents the solute region that includes the membrane and protein regions, while $\Omega_{s}$ represents the solvent and channel regions. $\Omega_{m}$ and $\Omega_{s}$ are separated by a boundary $\Gamma$. A two-dimensional illustration of the system in a simulation box is shown in Fig. \ref{Computationaldomain}.

Consider the following PNP equations \cite{lu2010FEM,Tu2013Aparallelfinite}:
\begin{equation*}
	\label{2.1-t}
	\left\{
	\begin{aligned} 
		&\frac{\partial p^i}{\partial t}-\nabla \cdot D^{i}\left(\nabla p^{i}+\beta q^{i} p^{i} \nabla \phi\right)=0, \text { in } \Omega_{s}, i=1, \ldots, N, \\  
		&-\nabla(\epsilon \nabla \phi)-\lambda \sum_{i=1}^{K} q^{i} p^{i}=\rho^{f}(x), \text { in } \Omega. \\ 
	\end{aligned}
	\right.
\end{equation*}
Here, $\phi$ and $p^i$ represent the electrostatic potential and the concentration of the $i$-th ionic species, respectively and $D^i$ denotes the diffusion coefficient of the $i$-th ionic species. The constant $\beta=k_BT$ is the inverse Boltzmann energy, where $K_B$ stands for the Boltzmann constant and $T$ is the absolute temperature. 
The characteristic function $\lambda$ equals 1 in $\Omega_{m}$ and 0 in $\Omega_{s}$, $\rho^f(x)=\sum_jq_j\delta\left(x-x_j\right)$ is an ensemble of fixed charges $q_j$ located at $x_j$. The piecewise dielectric constant $\epsilon(x)=\begin{cases}\epsilon_0\epsilon_m,    \text { in }\Omega_m,    \\\epsilon_0\epsilon_s,    \text { in }\Omega_s,    \end{cases}$
where $ \epsilon_0$ is the dielectric constant of vacuum. 


In this paper, since we focus on the space discretization, we only study the following steady-state PNP equations \cite{Tu2013Aparallelfinite, Wang2021}:
\begin{equation}
	\label{2.1}
	\left\{
	\begin{aligned} 
		&\frac{\partial p^i}{\partial t}-\nabla \cdot D^{i}\left(\nabla p^{i}+\beta q^{i} p^{i} \nabla \phi\right)=0, \text { in } \Omega_{s}, i=1, \ldots, N, \\  
		&-\nabla(\epsilon \nabla \phi)-\lambda \sum_{i=1}^{K} q^{i} p^{i}=\rho^{f}(x), \text { in } \Omega. \\ 
	\end{aligned}
	\right.
\end{equation}
The internal interface and boundary conditions are as follow:
\begin{equation}\label{2.2}
\begin{gathered}
\left[\phi\right]=0, \quad \left[\epsilon \frac{\partial \phi
}{\partial \mathbf{n}}\right]=0, \quad \text { on } \Gamma, \\
\phi=\phi_{D}, \quad \text { on } \Gamma_{D,\phi}, \quad\quad
p^i=p^{i}_D, \quad \text { on } \Gamma_{D,p^i}, \quad i=1, \ldots, N, \\
\frac{\partial \phi}{\partial n}=0, \quad \text { on } \Gamma_{N,\phi}, \quad\quad \mathbf{J_i}
 \cdot {\mathbf{n}}=0, \quad \text { on } \Gamma_{N,p^i}, \quad i=1, \ldots, N,
\end{gathered}
\end{equation}
where $\mathbf{J_i}=-D^i\left(\nabla p^i + \beta q^ip^i\nabla\phi\right)$, ${\mathbf{n}}$ represents the unit normal vector on the interface $\Gamma$, pointing from $\Omega_{m}$ to $\Omega_{s}$, $\Gamma_{D,\phi}$, $\Gamma_{N,\phi}$, $\Gamma_{D,p^i}$ and $\Gamma_{N, p^i}$ represent the Dirichlet boundary and the Neumann boundary for the electrostatic potential and ionic concentrations, respectively.

To deal with singular permanent charges in  (\ref{2.1}), the electrostatic potential $\phi$ in the Poisson equation is decomposed into three componets \cite{lu2010FEM}:
\begin{equation}\label{regularization}
\phi=\phi_{r}+\phi_s+\phi_m, 
\end{equation}
where $\phi_s=\frac1{4\pi\epsilon_ m}\sum_j\frac{q_j}{\left|x-x_j\right|},$
$\phi_m$ is the solution to the Laplace equation:
$$\begin{cases}-\Delta\phi_m=0,\quad \text{in} ~\Omega_m,\\\phi_m=-\phi_s, \text{ on } ~\partial \Omega_m.\end{cases}$$
The PNP equations consisting of the regularized Poisson equation is as follows: 

\begin{equation}
	\label{2.3}
	\left\{
	\begin{aligned} 
		&\nabla \cdot D^{i}\left(\nabla p^{i}+\beta q^{i} p^{i} \nabla \phi_r\right)=0, \text { in } \Omega_{s}, i=1, \ldots, N, \\  
		&-\nabla(\epsilon \nabla \phi_r)-\lambda \sum_{i=1}^{K} q^{i} p^{i}=0, \text { in } \Omega . \\ 
	\end{aligned}
	\right.
\end{equation}
The interface condition and the boundary conditions for the PNP equations (\ref{2.3}) are as follows:

\begin{equation}
\label{bc}
\begin{gathered}
\left[\phi_{r}\right] = 0, \quad \left[ \epsilon \frac {\partial \phi_r} {\partial \mathbf{n}}\right]=\epsilon_m \frac {\partial (\phi_s+\phi_m)}{\partial \mathbf{n}}  \quad \text{on } \Gamma, \\
\phi_r = \phi_{D}, \quad \text{on } \Gamma_{D, \phi}, \quad\quad
p^i = p^{i, D}, \quad \text{on }  \Gamma_{D, p^i}, \quad i = 1, \ldots, N, \\
\frac{\partial \phi_r}{\partial \textbf{n}} = 0, \quad \text{on } \Gamma_{D, \phi}, \quad\quad
\mathbf{J_i} \cdot \mathbf{n} = 0, \quad \text{on } \Gamma_{N, p^i}, \quad i = 1, \ldots, N.
\end{gathered}
\end{equation}
For simplicity, in the following we still use $\phi$ to denote the potential $\phi_r$. 

Define the following Sobolev spaces:
$$\begin{gathered}
H_{\phi}^1(\Omega)= \left\{v\mid v\in H^1(\Omega), v\mid_{\Gamma_{D,\phi}}=\phi_D\right\}, 
H_{p^i}^1(\Omega_s)= \left\{v\mid v\in H^1(\Omega_s), v\mid_{\Gamma_{D, p^i}}=p^{i}_D\right\},~~i=1,2,\cdots,N. \\
\end{gathered}$$
The weak formulation of (\ref{2.3}) is as follows: find $\phi\in H_\phi^1(\Omega),$ and $ p^i\in H_{p^i}^1(\Omega_s)$, satisfying:

\begin{equation}
	\label{2.4}
	\left\{
	\begin{aligned} 
		&(D^i\triangledown p^i, \triangledown v)+(\beta q^ip^i\triangledown\phi, \triangledown v)=0, ~~\forall v\in H_0^1(\Omega_s),    \\  
		&\left(\epsilon\triangledown\phi,    \triangledown w\right)-(\lambda\sum_{i=1}^Kq^ip^i,    w)+\left\langle\varepsilon_{m}\frac{\partial (\phi_{s}+\phi_m)} {\partial \mathbf{n}} , w\right\rangle_{\Gamma}=0, \forall w\in H_0^1(\Omega). \\ 
	\end{aligned}
	\right.
\end{equation}

As introduced in the Introduction, FEM has show the efficiency in solving PNP ion channel models. Since our primary focus is to improve the computational efficiency of the Gummel iteration combined with FEM discretization, we adopt the classical finite element discretization in this paper. For other improved FEMs, the procedure of the Gummel iteration remains nearly identical, so the proposed algorithm is also applicable to these other improved finite element discretizations.

Suppose $T^h(\Omega)$ is a mesh of size $h$ on the domain $\Omega$. Define the following linear finite element spaces:
$$\begin{gathered}
V_{h}  = \left\{v_{h}\left|v_{h} \in H^{1}(\Omega), v_{h}\right|_{e} \in P^{1}(e), \forall \mathrm{e} \in T^{h}(\Omega),    \left.v_{h}\right|_{\Gamma_{D, \phi}}  =  \phi_{D}\right\}, \\
V_{h}^{0}  =  \left\{v_{h}\left|v_{h} \in H^{1}(\Omega), v_{h}\right|_{e} \in P^{1}(e), \forall \mathrm{e} \in T^{h}(\Omega),    \left.v_{h}\right|_{\Gamma_{D, \phi}}  =  0\right\}, \\
S_{h}  =  \left\{w_{h}\left|w_{h} \in H^{1}(\Omega_s), w_{h}\right|_{e} \in P^{1}(e), \forall \mathrm{e} \in T^{h}(\Omega),    \left.w_{h}\right|_{\Gamma_{D, p^i}}  =  p^{i}_D\right\}, \\
S_{h}^{0}  =  \left\{w_{h}\left|w_{h} \in H^{1}(\Omega_s), w_{h}\right|_{e} \in P^{1}(e), \forall \mathrm{e} \in T^{h}(\Omega),    \left.w_{h}\right|_{\Gamma_{D, p^i}}  =  0\right\} .
\end{gathered}$$
The standard finite element approximation corresponding to equation (\ref{2.4}) is defined as:
finding $\phi_{h}\in V_{h}$ and $ {p}_{h}^{i}\in S_{h}, {i}={1}, ..., {N} ,$ such that:

\begin{equation}
	\label{2.5}
	\left\{
	\begin{aligned} 
		&(D^i\nabla p_h^i,    \nabla v_h)+(\beta q^ip_h^i\nabla\phi_{h},    \nabla v_h)=0, \forall v_h\in V_h^0,    \\  
		&(\epsilon\nabla\phi_h,    \triangledown w_h)-(\lambda\sum_{i=1}^Kq^ip_h^i,    w_h)+\left\langle\varepsilon_{m}\frac{\partial (\phi_{s}+\phi_m)} {\partial \mathbf{n}} , w_h\right\rangle_{\Gamma}=0, \forall w_h\in S_h^0.\\ 
	\end{aligned}
	\right.
\end{equation}

\setcounter{lemma}{0}
\setcounter{theorem}{0}
\setcounter{corollary}{0}
\section{Selection of relaxation parameters for Gummel iteration based on RFR Method}

In this section, we first introduce the FEM with Gummel iteration for the PNP problem \eqref{2.3}-\eqref{bc}, and the entire solution procedure for solving the system. After that, we introduce the RFR method and propose three algorithms to generate the training set for the RFR method. Finally, we present an novel iterative algorithm based on RFR method for this PNPi model.



\subsection{The Gummel iteration}\setcounter{equation}{0}

Due to the nonlinearity and strong coupling characteristics of PNP equations in ion channels, the decoupling strategies are commonly employed in practice computation. The Gummel iteration method, as an efficient and simple decoupling approach, has been successful applied to simulate several ion channels \cite{Tu2013Aparallelfinite, Zhang2020, Wang2021}. The Gummel iteration was originally proposed by Gummel \cite{gu64} to solve the drift-diffusion equation in semiconductor devices. Chen et al. provided a general Gummel iteration for the PNP equations \cite{chen2013surface}. This iteration is a simple step-by-step formula, which solves one type of equation (for example, the Poisson equations) first and then substitutes the obtained solution into another type equation to decouple the system. 



We apply the Gummel iteration \cite{chen2013surface} to the finite element formulation \eqref{2.5} to present the FEM with Gummel iteration algorithm for the PNP problem (\ref{2.3})-(\ref{bc}) , as shown in Algorithm \ref{Gummel}.
\begin{algorithm}  
\caption{FEM with Gummel iteration for PNP equations}  
\label{Gummel}

{\bf Step 1.} Give the initial error tolerance $tol$, set the maximum number of iterations, and initialize the relaxation parameter $\alpha$. Set initial iterative step counter $k=0$ and initial concentration $p_h^{i,k} =0$\\

{\bf Step 2.} Solve the decoupled PNP.\\
\While{$k\le K$ } 
{ Find $\phi_h^{k+1}\in S_h$ and $p_h^{k+1}\in V_h$, such that
$(\epsilon\nabla\phi_{h}^{k+1},    \triangledown w_h)-(\lambda\sum_{i=1}^Nq^i p_h^{i,    k},    w_h)+\left\langle\varepsilon_{m}\left(\frac{\partial \phi_{s}} {\partial n^m}+\frac{\partial \phi_{h}} {\partial n^m}\right), w_h\right\rangle_{\Gamma}=0, \forall w_h\in S_h^0\;\label{3.2}$,\\

	$(D^i\nabla p_h^{i,    k+1},    \nabla v_h)+(\frac e{K_BT}q^i p_h^{i,    k+1}\nabla\phi_{h}^{k+1},    \nabla v_h)=0, \forall v_h\in V_h^0\;\label{3.1}$.
	
	\If{$\left | \left | \phi_h^{k+1} - \phi_h^{k}\right | \right | + \left | \left | p_h^{i,    k+1} - p_h^{i,    k}\right | \right | < tol$}{
		\textbf{break}\;  
	}
	\Else{
		$\phi_{h}^{k+1} \gets (1-\alpha) \phi_{h}^{k} + \alpha \phi_{    h}^{k+1}$\;
		$p_h^{i,    k+1} \gets (1-\alpha) p_h^{i,    k} + \alpha p_h^{i,    k+1}$\;
	}
	$k \gets k + 1$\;
}
{\bf Step 3.} Output the electrostatic potential $\phi_{h}^{k+1}$ and concentration $ p_h^{i, k+1},~i=1,2,\cdots,N$.

\end{algorithm} 

To ensure the convergence of the Gummel iteration,
the relaxation parameter $\alpha$ (where $0 <\alpha< 1$) is always required \cite{lu2010FEM,    Tu2013Aparallelfinite}, see it in Step 2 of Algorithm \ref{Gummel}. The parameter $\alpha$ is a value between $0$ and $1$. For practical PNP problems, the value of the relaxation parameter $\alpha$ has a significant influence on the efficiency of the iterative method.
This is because, first, the convergence efficiency of Algorithm \ref{Gummel} is sensitive to the value of $\alpha$, which can be seen from Table... Second, the solution of practical PNP problems usually requires fine computation within a large-scale grid domain with extremely small grid sizes, which leads to very high numerical computation costs. If the relaxation parameters are improperly selected, it will directly increase the number of iterations, and each additional iteration will further consume a large amount of computational resources on top of the already high computational costs.
  
Traditionally, the selection of relaxation parameters mostly relies on experience. However, parameter values selected based on manual experience are inappropriate or not stable for practical PNP equations. As shown the numerical experiments in Section 4, the manually selected parameter is not an optimal one and sometime leads to divergence. In next section, we shall propose a parameter selection method based on a machine learning.

\begin{figure}[h]
	\centering
	\includegraphics[height=5cm, width=16cm]{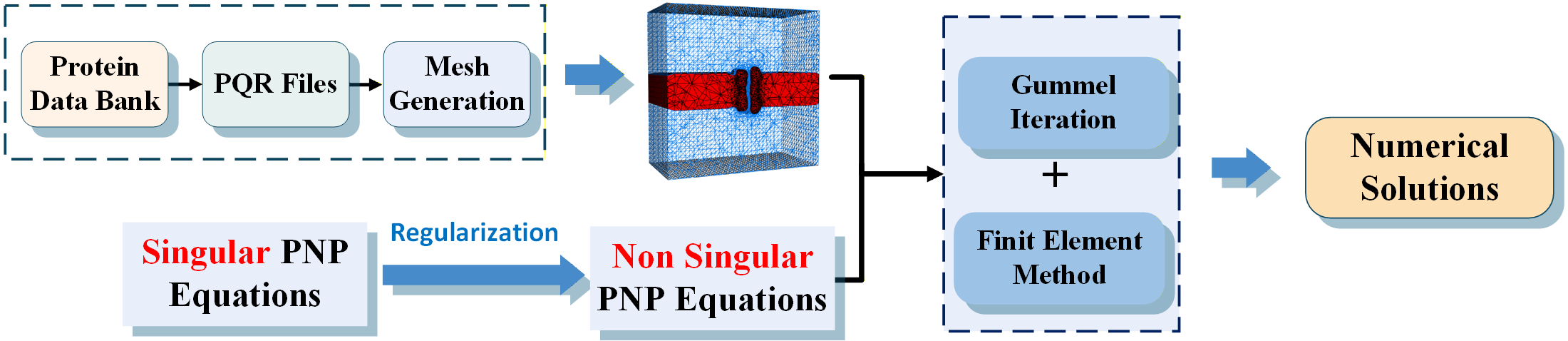}
	\caption{Solution procedure for the PNP equations in ion channels.}
	\label{FEMsolvePNP}
\end{figure}

Based on Algorithm \ref{Gummel}, we can present the whole solution procedure for solving PNP equations \eqref{2.1} satisfying conditions \eqref{2.2}, which is demonstrated in Fig. \ref{FEMsolvePNP}. The whole solution procedure includes three main steps: First, find the PQR structure files of the target ion channel in the protein bank and use the mesh generation platform/software TMSmesh and TetGen \cite{Tu2013Aparallelfinite} to generate the mesh for the ion channel system. Second, to address the singularity in the governing equations (\ref{2.1})–(\ref{2.2}), the regularization method given by equation (\ref{regularization}) is applied to transform them into the forms  (\ref{2.3})–(\ref{bc}). Finally, the nonsingular PNP problem  (\ref{2.3})–(\ref{bc}) is solved using Algorithm \ref{Gummel} to get the numerical solution.

\subsection{ RFR model and training datasets generation}

RFR is a machine learning method based on ensemble learning. It can effectively suppress noise and overfitting, thereby capturing the complex nonlinear mapping between features and the target variable. RFR has been successfully applied to parameter selection for the Picard iteration in reservoir models \cite{Silva2021Machinelearning}. In this paper, we apply RFR to parameter selection for the Gummel iteration in biological ion channel simulations. We also employed another machine learning model, namely Gaussian Process Regression. However, experimental results show that its predictive performance is inferior to that of RFR. Therefore, only the RFR model is introduced in this work.

\begin{figure}[H]
	\centering
	\includegraphics[height=8cm,    width=12cm]{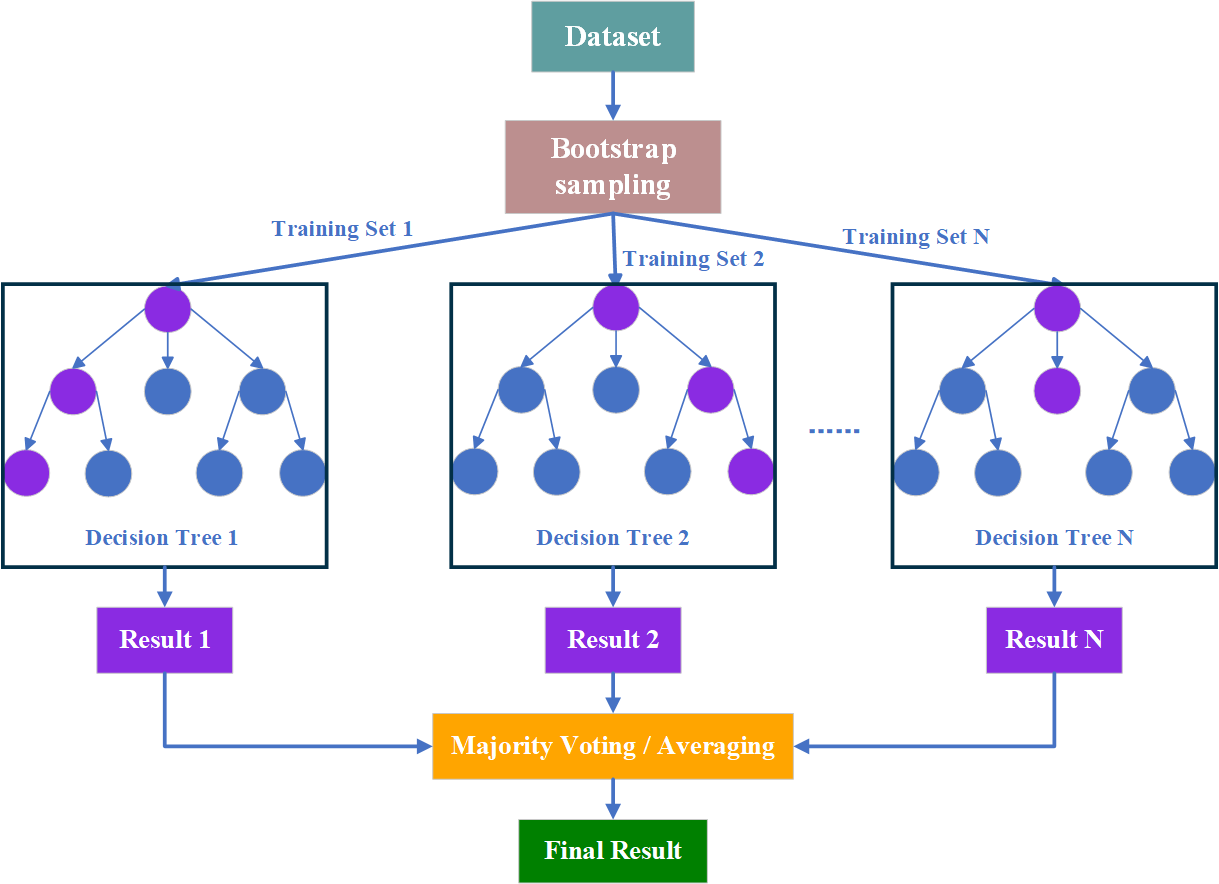}
	\caption{Random Forest Regression model schematic}
	\label{RFR_tree}
\end{figure}

The RFR was initially proposed by Breiman in 2001\cite{breiman2001random} and serves as a typical representative of the Bagging algorithm within ensemble learning. Its “randomness” is manifested in two aspects. Firstly, at the sample level, sub-sample sets are randomly drawn with replacement from the original dataset via Bootstrap sampling. Secondly, at the feature level, when splitting each node of a decision tree, a subset of the features is randomly selected to determine the optimal splitting point. The term “forest” refers to the integration of multiple decision trees constructed based on the CART algorithm. The specific process is illustrated in Fig. \ref{RFR_tree}, for which the dataset must be pre-generated when employing the RFR model.

To construct the training set for the RFR model in ion channel simulation, it is necessary to incorporate physical, geometrical and numerical quantities that have a crucial impact on the behavior of ion channels as input features (a detailed introduction is given in Table \ref{tab:input_parameters} in Section 4). Some of the quantities are easy to obtain, such as voltage, ion concentration, and charge density on the membrane. However, two important quantities are challenging to generate. One is the information related to the numerical solution, such as the number of iterations. The other is the geometrical shape features. Next, we describe these two challenges and present the algorithms used to address them.

Due to the complexity of numerical simulations for ion channels, obtaining numerical solutions for channels with complex structures is often difficult. In existing work on the numerical simulation of ion channels, only a very small number of channels with relatively simple geometries can yield numerical solutions.
For complex large channels, numerical calculations struggle to converge, making numerical solutions extremely difficult to obtain. This results in a very sparse training dataset of numerical solutions available for the RFR model.

\textcolor{black}{To address this issue, we draw on the research approach of \cite{Mohammad2010}. In that study, FhuA-like proteins were redesigned to construct large nanopores suitable for single-molecule sensing. Since FhuA possesses a stable $\beta$-barrel transmembrane scaffold and a relatively large internal pore, and its structure can be readily simplified through residue deletion, it is well suited as a basic template for constructing artificial ion channels. Based on this idea, we propose an artificial channel generation algorithm, as detailed in Algorithm \ref{artificial channel generation}. This algorithm can generate a class of simplified and numerically tractable artificial ion channels based on FhuA-like proteins, thereby providing sufficient numerical solution data.}

Algorithm \ref{artificial channel generation} aims to modify the PQR files \textcolor{black}{of the FhuA protein} and generate new PQR files for artificial channels, which forms the foundation for mesh generation in ion channel simulations (cf. Fig. \ref{FEMsolvePNP}). The PQR files for natural channels are molecular structure files widely used in computational biochemistry and structural biology. These files are usually generated by third-party software packages such as PDB2PQR and APBS. The specific format of the PQR file is presented in Table \ref{tab:pqr_format}.

 \begin{table}
\centering
\caption{Format and variables of the PQR file}
\label{tab:pqr_format}
\begin{tabular}{cll}
\toprule
\textbf{Columns} & \textbf{Variable} & \textbf{Definition} \\
\midrule
1-6   & Record type & Entry type, should be "ATOM" or "HETATM" \\
7-11  & Atom serial number & Unique index of the atom \\
12-16 & Atom name & Element symbol or atom identifier (e.g., CA, N, C, O) \\
17-20 & Residue name & Name of the residue (e.g., ALA, LYS) \\
21-22 & Chain identifier & Chain ID of the residue \\
23-26 & Residue number & Sequence number of the residue \\
31-38 & X coordinate & Atomic coordinate in X (Å) \\
39-46 & Y coordinate & Atomic coordinate in Y (Å) \\
47-54 & Z coordinate & Atomic coordinate in Z (Å) \\
55-62 & Charge (q) & Partial charge of the atom\\
63-69 & Radius (r) & Atomic radius (Å) \\
\bottomrule
\end{tabular}
\end{table}

\textcolor{black}{The core steps of Algorithm \ref{artificial channel generation} are as follows. First, the FhuA protein structure is used as the initial template, and the set of residues to be removed from the PQR file, together with their spatial range, is determined according to the residue deletion strategy proposed in \cite{Mohammad2010}. The corresponding residues in the original FhuA structure are then removed to obtain an artificial channel model with a more open pore cavity and a simplified structure. During this process, the algorithm preserves, as much as possible, the original transmembrane scaffold, electrostatic features, and geometric constraints of FhuA, while reducing the structural complexity of the natural channel system. This provides biophysically reasonable artificial ion channel structures for subsequent numerical simulations. Then, the visualization software VMD (Visual Molecular Dynamics) is used to precisely locate the remove residue set for constructing artificially generated channels. VMD is a widely used software tool for molecular structure visualization and analysis \cite{HUMPHREY199633}, supporting various file formats such as PQR and PDB, and enabling three-dimensional spatial localization, filtering, and export of atomic subsets.}

\begin{algorithm}[H]  
	\caption{Artificial channel generation base on FhuA protein}  
	\label{artificial channel generation} 
	\textbf{Input:} PQR file of the wild-type FhuA protein\\
	\textbf{Output:} New PQR file of the artificial channel
	
	\textbf{Step 1: Original Channel Data Input} 
    \begin{itemize}
	    \item Read and import the PQR file of the wild-type FhuA protein. For each atom, extract its three-dimensional coordinates along with its residue information (residue name, residue number, chain identifier) and physical parameters (charge $q_i$, radius $r_i$). Store the complete data as a set of tuples $\{(x_i, y_i, z_i, \text{res}_i, q_i, r_i)\}_{i=1}^{N_a}$, where $\text{res}_i$ denotes the parent residue of atom $i$ and $N_a$ is the total number of atoms.
	\end{itemize}
	
      \textbf{Step 2: Generation of Atom Subsets ${P}_{\text{new}}$ Based on Residue Deletion Strategy}
      \begin{itemize}
          \item Based on the deletion strategy reported in \cite{Mohammad2010} for FhuA, construct the set of residues to be removed as $\mathcal{R}_{\text{remove}}$. This set includes the plug domain and four extracellular loops (L3, L4, L5, and L11), with specific residue ranges as defined in \cite{Mohammad2010}.
          \item        
          For each residue \( R_j \in \mathcal{R} \), use the visualization software VMD \cite{HUMPHREY199633} to locate and extract its corresponding atom subset \( P_j \), and assisted selection through visual interface to achieve more accurate localization.
        \item The atoms to be removed are then the union of atom subsets corresponding to these residues:
     $${P}_{\text{remove}} = \bigcup_{R_j \in \mathcal{R}_{\text{remove}}} {P}_j.$$
     \item The updated atom set after removal is given by ${P}_{\text{new}} = {P}    \setminus {P}_{\text{remove}}$, where P is the set of all atoms.
    \end{itemize}
    
	\textbf{Step 3: Artificial Channel Data Output} 
    \begin{itemize}
	    \item For each atom in the updated set ${P}_{\text{new}}$, retrieve its corresponding charge $q_i$ and radius $r_i$ from the original PQR file.
	    \item Write the atomic coordinates along with their charges and radii into a new PQR file, following the standard PQR column format (see Table \ref{tab:pqr_format} for details).
	    \item Output the resulting PQR file for subsequent numerical simulations.
    \end{itemize}
\label{Artificial}    
\end{algorithm}
\vspace{1em}

\begin{figure}[!htbp]	
\begin{subfigure}[htb]{0.5\linewidth}		
\centering		
\includegraphics[height=5cm,width=4cm]{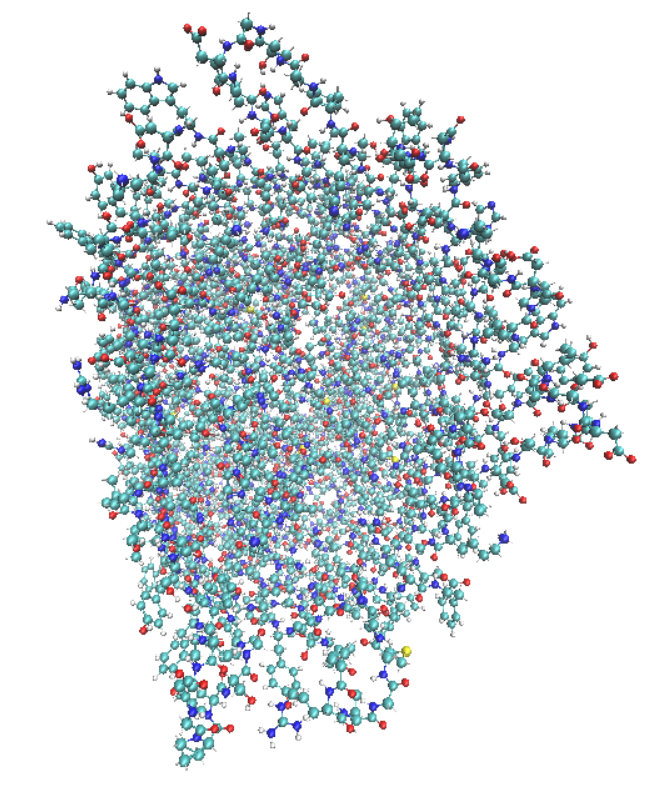}		
\caption{Structure of the original FhuA protein}  
\label{fig:FhuA_org}	
\end{subfigure}	
\begin{subfigure}[htb]{0.5\linewidth}		
\centering		
\includegraphics[height=5cm,width=4.8cm]{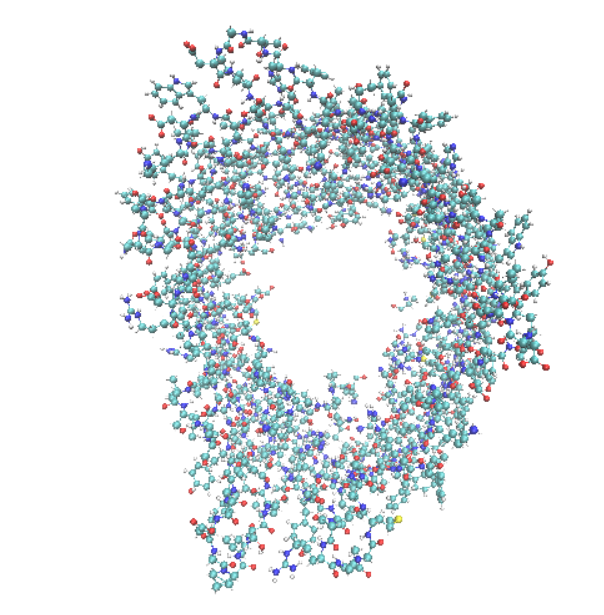}		
\caption{Structure of the FhuA protein after Algorithm \ref{artificial channel generation}}  
 \label{fig:FhuA_Alg}	
 \end{subfigure}
	
 \caption{ The structures of the original FhuA protein and the FhuA protein after Algorithm \ref{artificial channel generation}}	
 \label{fig:FhuA}
\end{figure}

\textcolor{black}{Fig. \ref{fig:FhuA} shows the atomic structure of the original FhuA protein and that after removing the residue set $\mathcal{R}_{\text{remove}}$ (including the cork domain and four extracellular loops L3, L4, L5, and L11) using Algorithm \ref{artificial channel generation}. A clear pore can be observed in the center of the FhuA protein after deletion, allowing ions in the solution to move freely through it.}

  Another challenge in the practical construction of training sets is how to extraction geometric shape features of channels effectively. The initial geometric structure data of ion channels were obtained from the RCSB Protein Data Bank (https://www.rcsb.org) in the PQR file format. 
The amount of PQR data is enormous. Directly introducing these atomic data into the RFR model may introduce a large number of redundant features, thereby increasing the computational burden. Inspired by the point cloud algorithm \cite{pfister2000surfels}, we propose two channel geometric feature extraction algorithms--Algorithms \ref{for Simple Channel} and \ref{for Complex Channel}. Algorithm \ref{for Simple Channel} is suitable for channels with relatively simple shapes. Built upon Algorithm \ref{for Simple Channel}, Algorithm \ref{for Complex Channel} is designed for channels with more complex geometries, yet it is less computationally efficient than Algorithm \ref{for Simple Channel}.

The main procedure of Algorithm \ref{for Simple Channel} is as follows: First, the voxel grid downsampling method \cite{rusu20113d} is applied to process the original PQR file. By partitioning the data into regular voxel grids and aggregating the points within each voxel, the data volume is  substantially reduced while the original geometric structure is preserved. Second, to extract effective geometric features, the K-D tree algorithm \cite{bentley1975multidimensional} is adopted to extract the channel surface region. Finally, the Farthest Point Sampling (FPS) algorithm \cite{eldar1997farthest} is used to progressively select a representative point set within this surface region, serving as the key geometric feature region for describing the shape of the channel protein.

\begin{algorithm}[!htbp]
	\caption{Geometric Feature Extraction Algorithm I }
	\label{for Simple Channel}
	
	\textbf{Input}: PQR file containing relevant channel data and \textcolor{black}{initial voxel size $r$}.
	
	\textbf{Output}: Geometric feature $G_{train}$ of the channel.
	
	\textbf{Step 1: Coordinate Data Extraction}
	\begin{itemize}
		\item Read the input PQR file and extract three-dimensional coordinate data. Convert the extracted coordinate data into point cloud objects using encapsulated functions from the Open3D library in Python, denoted as $P = \{\mathbf{p}_i\}_{i=1}^N = \{(x_i, y_i, z_i)\}_{i=1}^N$, where $N$ is the total number of atoms.
	\end{itemize}	
	
	\textbf{Step 2: Voxel Downsampling}
	\begin{itemize}
		\item Apply voxel downsampling \cite{rusu20113d} to the point cloud $P$ using a voxel size $r$. The 3D space is partitioned into a regular grid of cubic voxels with side length $r$. For the $k$-th non-empty voxel, let the set of points falling within it be denoted as:
		\[
		\mathcal{P}_k = \{ \mathbf{p} \in P \mid \mathbf{p} \text{ lies in the } k\text{-th voxel} \}.
		\]
		The representative point for this voxel is defined as the centroid of the points in $\mathcal{P}_k$:
		\[
		\mathbf{v}_k = \frac{1}{|\mathcal{P}_k|} \sum_{\mathbf{p} \in \mathcal{P}_k} \mathbf{p},
		\]
		where $|\mathcal{P}_k|$ denotes the number of points in the $k$-th voxel. The downsampled point cloud is then given by:
		$V = \{ \mathbf{v}_1, \mathbf{v}_2, \ldots, \mathbf{v}_M \},$
		where $M$ is the total number of non-empty voxels.  
	\end{itemize}
	
	\textbf{Step 3: Neighborhood Density Calculation and Low-Density Point Set Selection}
	\begin{itemize}
		\item For each point $\mathbf{v}_i \in V$, a K-D tree \cite{bentley1975multidimensional} is employed to accelerate the search for its $K$ nearest neighbors. The local neighborhood density $\rho_i$ for point $\mathbf{v}_i$ is then defined as:
		\[
		\rho_i = \frac{K}{\text{AvgDist}_i},
		\]
		where $\text{AvgDist}_i$ is the average distance to its $K$ nearest neighbors.
		\item Set a density threshold $\tau$. Identify all the points with a density below this threshold to form the low-density point set:
		\[
		V_{ld} = \{ \mathbf{v}_i \in V \mid \rho_i < \tau \}.
		\]
	\end{itemize}
	
	\textbf{Step 4: Key Geometric Feature Screening (continued)}
	\begin{itemize}
		\item \textcolor{black}{The FPS algorithm \cite{eldar1997farthest} is applied to \( V_{ld} \) to screen and extract the key geometric features of the channel. First, randomly select an initial point \( s_1 \) from \( V_{id} \) and set \( S = \{s_1\} \). Then, iteratively perform the following: at step \( t \), select a point \( s_t \) from \( V_{id} \setminus S \) that satisfies:
		\[
		s_t = \arg\max_{p \in V_{id} \setminus S} \ \min_{q \in S} \|p - q\|_2
		\]	
		where \( p \) is a currently unselected point, \( q \) is an arbitrary point in the already selected set \( S \), and \( \|p - q\|_2 \) denotes the Euclidean distance between the two points. 
        Repeat this process until \( |S| = g \). The final geometric feature set of the channel is obtained as:}
		
		\[
		G_g = S = \{g_i\}_{i=1}^g = \{(x_i, y_i, z_i)\}_{i=1}^g
		\]
		
		where \( G_g \subseteq V_{id} \), and \( g \) is a preset number of sampled points.
	\end{itemize}
    
	\begin{itemize}
		\item Output the geometric feature $G_{train}$ . 
	\end{itemize}
\end{algorithm}



For channels with complex structures, such as those containing a large number of atoms and highly irregular surfaces, directly applying Algorithm \ref{for Simple Channel} makes it difficult to accurately capture local features such as depressions and protrusions. To address this issue, we propose Algorithm \ref{for Complex Channel}, which is built upon Algorithm \ref{for Simple Channel}. The core improvement lies in the introduction of the the Random Sample Consensus (RANSAC) algorithm \cite{fischler1981random}: First, this algorithm is used to segment the complex channel protein into multiple local fragments with gentle geometric variations. Then, each fragment is independently subjected to geometric fitting to capture its local details. Finally, Steps 2–5 of Algorithm \ref{for Simple Channel} are applied to extract and integrate the geometric features of each fragment. By decomposing the global complex problem into a series of local simple problems, the new algorithm significantly enhances its adaptability to complex geometric structures.

It is worth noting that Algorithm \ref{for Complex Channel} is designed for channels with complex geometries (test channels VDAC and KcsA), while Algorithm \ref{for Simple Channel} applies to channels with simple geometries (training channel GA and artificially generated channels).
\vspace{1em} 

\begin{algorithm}[H]
	\caption{Geometric Feature Extraction Algorithm II}
	\label{for Complex Channel}
	
	\textbf{Input}: PQR file containing 3D coordinates of ion channel atoms.\\
	\textbf{Output}: Geometric feature vector $G_{test}$ representing the channel structure.
	
	\textbf{Step 1: Data Acquisition}
	\begin{itemize}
		\item Read the input PQR file and extract the atomic coordinate data $(x, y, z)$ for each atom in the ion channel structure.
	\end{itemize}

	\textbf{Step 2: Point Cloud Construction}
	\begin{itemize}
		\item Convert the extracted coordinate data into a point cloud object $P$ using the encapsulated functions from the Open3D library in Python.
	\end{itemize}
	
	\textbf{Step 3: Plane Segmentation via RANSAC}
	\begin{itemize}
		\item Apply the RANSAC algorithm \cite{fischler1981random} to segment the point cloud $P$ into planar components.
		\item Identify and store all planar regions representing the channel surface into a plane set $P_{\text{planes}} = \{P_1, P_2, \ldots, P_J\}$, where $P_j, j=1,2,...,J$ is a subset of points belonging to a distinct plane.
	\end{itemize}
	
	\textbf{Step 4: Geometric Feature Extraction}
	\begin{itemize}
		\item Process the plane set $P_{\text{planes}}$ using Steps 2-4 of Algorithm \ref{for Simple Channel} to obtain geometric feature vector $G_{test}$.
	\end{itemize}
	
\end{algorithm}
\vspace{1em}

\textcolor{black}{The training dataset for Algorithm \ref{RFR} consists of multidimensional feature vectors $(K_{\text{train}}, G_{\text{train}}, S)$. Here, $G_{\text{train}} = (G_{\text{nat}}, G_{\text{art}})$ denotes the geometric features of training channels extracted by Algorithm \ref{for Simple Channel}, including natural channels GA and artificially generated channels. The physicochemical parameter vector $K_{\text{train}} = {(V_i, C_i, \sigma_i)}_{i=1}^{N}$ is composed of boundary information voltage (V), initial ion concentration (C) and equation parameters (membrane charge density $\sigma$). $S = {(\alpha_i, \text{iters}i)}_{i=1}^{N}$ represents the numerical solution information, where $\text{iters}_i$ is the number of iterations obtained by substituting the relaxation parameter $\alpha_i$ into Algorithm \ref{Gummel}.} The performance of Algorithms \ref{for Simple Channel} and \ref{for Complex Channel} will be shown in the numerical experiments of Section 4.

\pagebreak\noindent
\subsection{An improved Gummel iteration based on RFR parameter selection}
This section presents Algorithm \ref{RFR}, an improved Gummel iteration based on the RFR parameter selection. The algorithm first trains the RFR model using datasets constructed by Algorithms \ref{artificial channel generation}, \ref{for Simple Channel}, and \ref{for Complex Channel}. Building on this, the trained RFR model is integrated with Algorithm \ref{Gummel} to accelerate the iterative solution of the PNP equations by optimizing the relaxation parameters. The key distinction from Algorithm \ref{Gummel} is that the proposed method employs RFR to automatically determine appropriate or near-optimal relaxation parameters, rather than relying on empirical manual selection, thus enhancing the convergence efficiency of the iterative process.

\begin{figure}[h]
\centering
\includegraphics[height=8.5cm,    width=15cm]{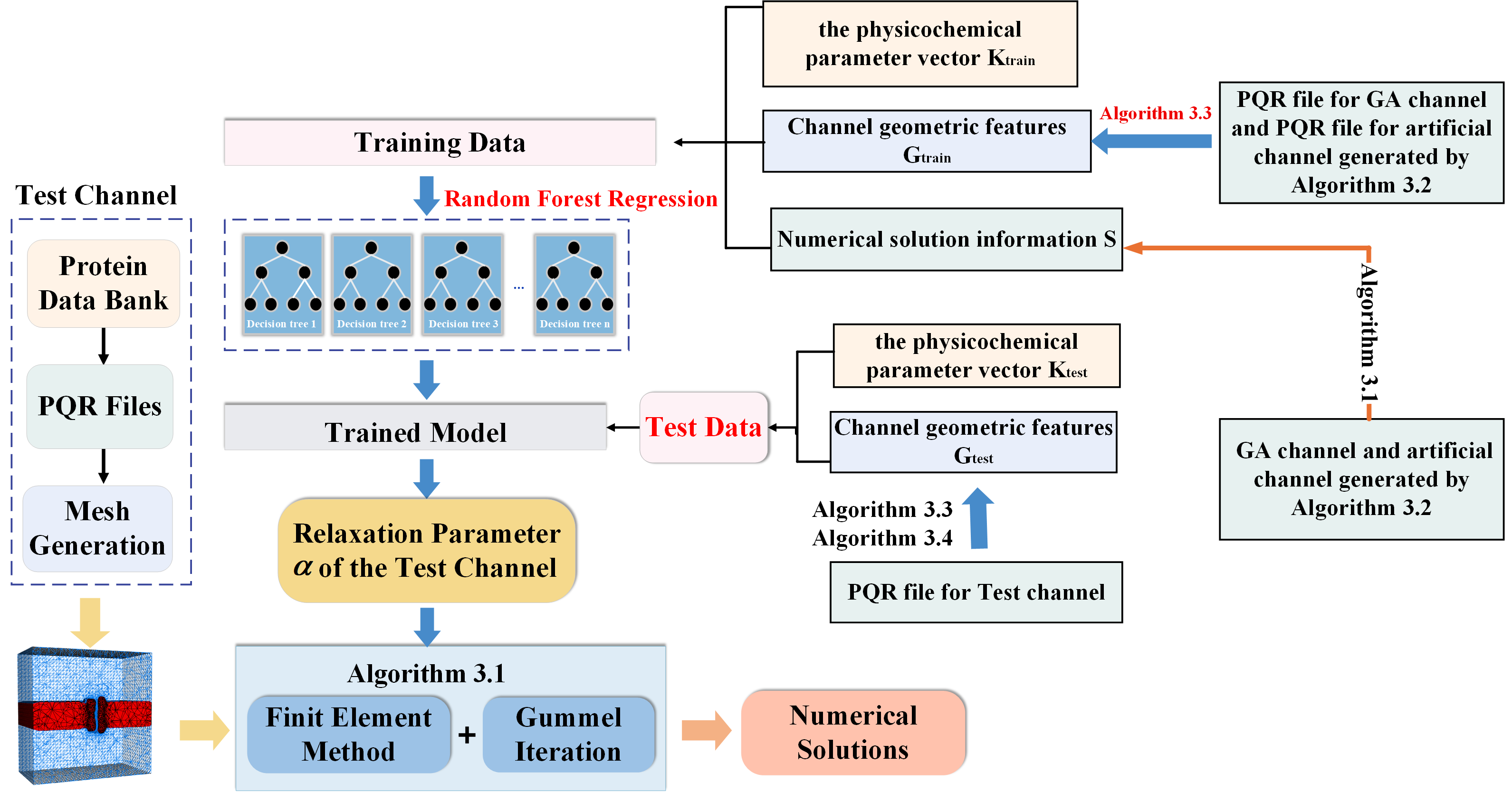}
\caption{Solution process of Algorithm \ref{RFR} for the PNP equations in ion channel simulations }
\label{RFR+FEM}
\end{figure}

\textcolor{black}{Fig. \ref{RFR+FEM} illustrates the overall workflow for solving the PNP equations in ion channel simulation based on Algorithm \ref{RFR}. First, Algorithm~\ref{artificial channel generation} is employed to construct artificial ion channels to augment the training dataset. Second, Algorithm~\ref{for Simple Channel} is used to extract the geometric features $G_{\text{train}} = (G_{\text{nat}}, G_{\text{art}})$ of the natural channel GA and the artificial channels. Third, a set of relaxation parameters $\alpha_i$ are respectively substituted into Algorithm~\ref{Gummel} to solve the PNP equations for the above two types of channels and obtain the corresponding iteration numbers $\text{iters}_i$, which together with $\alpha_i$ form the data pairs $S = \{(\alpha_i, \text{iters}_i)\}_{i=1}^{N}$. Then, these are combined with a set of physical parameters $K_{\text{train}} = \{(V_i, C_i, \sigma_i)\}_{i=1}^{N}$ (where $V_i$, $C_i$, and $\sigma_i$ denote the voltage, ion concentration, and membrane charge density, respectively) to construct the training dataset $(K_{\text{train}}, G_{\text{train}}, S)$ for training the RFR model. Finally, the trained model is applied to the geometric features $G_{\text{test}}$ and physical parameters $K_{\text{test}}$ of the test channels to adaptively select the relaxation parameters for the Gummel iteration, thereby achieving an efficient solution of the PNP equations in conjunction with the finite element method.}

\vspace{1em}

\vspace{1em}

\begin{algorithm}[H]
	\caption{An improved Gummel iteration based on RFR parameter selection
}
	\label{RFR}
	
	\textbf{Input}: A set of training physicochemical parameter vectors $K_{\text{train}} = \{(V_i, C_i, \sigma_i)\}_{i=1}^{N}$, a test physicochemical parameter vector $K_{\text{test}} = (V, C, \sigma)$, the atomic coordinate data from the training PQR files $Q_{\text{train}} = (Q_{\text{nat}}, Q_{\text{FhuA}})$, the atomic coordinate data from the test PQR file $Q_{\text{tst}}$, and a set of relaxation parameters $\alpha$.
	
	\textbf{Output}: Numerical solutions of the PNP equations for the ion channel.
	
	\textbf{Step 1: Training Set Generation}
	\begin{itemize}
		\item 
        Using Algorithm~\ref{artificial channel generation}, the atomic coordinate data $Q_{\text{FhuA}}$ is processed to generate a new PQR file $Q_{\text{art}}$ for the artificial channel.

       \item 
       $Q_{\text{nat}}$ and $Q_{\text{art}}$ are each processed using TMSmesh \cite{Tu2013Aparallelfinite} to generate their corresponding meshes. Then, based on the physicochemical parameter vector $K_{\text{train}} $ and the relaxation parameter $\alpha$, the corresponding number of iterations $\text{iters}$ is computed via Algorithm \ref{Gummel}, which together with $\alpha$ forms the data pair $S = (\alpha, \text{iters})$.

        \item 
        Using Algorithm \ref{for Simple Channel} or \ref{for Complex Channel}, geometric feature extraction is performed on $Q_{\text{nat}}$ and $Q_{\text{art}}$ respectively, yielding the geometric feature set $G_{\text{train}} = (G_{\text{nat}}, G_{\text{art}})$
		\item 
        In summary, the multi-dimensional feature vector $(K_{\text{train}}, G_{\text{train}}, S)$, composed of the physicochemical parameter vector $P$, the geometric features $G$, and the data pair $S$, is used as the training dataset.
        
	\end{itemize}
	
	\textbf{Step 2: RFR Modeling}
	\begin{itemize}
		\item Initialize the RFR model \cite{breiman2001random} by setting the number of trees as $T$, max features per node $m$, and the maximum depth of the trees as $depth$.
		\item Input the training dataset constructed in Step 1 into the RFR model and train the model to obtain a trained RFR model. 
	\end{itemize}
	
	\textbf{Step 3: Prediction of Relaxation Parameter}
	\begin{itemize}
		\item 
        Extract the geometric features $G_{\text{test}}$ of the test channel (using Algorithms~\ref{for Simple Channel}–\ref{for Complex Channel}) together with the physicochemical parameter vector $K_{\text{test}}$ to form the test dataset.
        
		\item Input the test dataset into the trained RFR model. The model will then predict and output the corresponding relaxation parameter $\alpha$.
	\end{itemize}
	
	\textbf{Step 4: Solving PNP equations}
    \begin{itemize}
   \item Apply the  relaxation parameter $\alpha$ obtained in Step 3 to Algorithm \ref{Gummel}, and output the numerical solutions $\phi_h$ and $p_h^i,~~i=1,2,\cdots,M$, and the number of iteration of Algorithm \ref{Gummel}.
    \end{itemize}
	
\end{algorithm}
\vspace{1em} 

\textbf{Remark:} $K_{\text{train}} = \{(V_i, C_i, \sigma_i)\}_{i=1}^{N}$, where $V_i$, $C_i$, and $\sigma_i$ denote the voltage, ion concentration, and membrane charge density, respectively; $Q_{\text{train}} = (Q_{\text{nat}}, Q_{\text{art}})$, where $Q_{\text{nat}}$ is the atomic coordinate data of the GA channel and $Q_{\text{art}}$ is the atomic coordinate data of the FhuA protein.

\section{Numerical example}\setcounter{equation}{0}

Consider the test channels immersed in a 1:1 KCl solution. The channel contains two ionic species: monovalent potassium ions ($\mathrm{K}^+$) and monovalent chloride ions ($\mathrm{Cl}^-$). Let $p^1$ and $p^2$ denote the concentrations of $\mathrm{K}^+$ and $\mathrm{Cl}^-$ in the bulk solution, respectively. The ion transport process is described by the PNP equations (\ref{2.3}) together with their boundary and interface conditions (\ref{bc}). 

Here, $\phi$ is the electrostatic potential and the dielectric coefficient $\epsilon(x)=\left\{\begin{array}{l}{2\epsilon_{0}, x\in\Omega_{m}, }\\{80\epsilon_{0}, x\in\Omega_{s}.}\\\end{array}\right.$ Some of the main parameters in (\ref{2.3}) are reported in Table \ref{tab:ionparameters}.

\begin{table}[h]	
\centering	
\caption{The parameters of (\ref{2.1})-(\ref{2.1})}	
\label{tab:ionparameters}	
\begin{tabular}{c c p{6cm}}		
\hline		
\textbf{Symbol} & \textbf{Value} & \textbf{Description} \\		
\hline		
$D^1$ & $2.0561 \times 10^{-9} \, \text{m}^2/\text{s}$ & Diffusion coefficient ($\mathrm{K}^+$) \\		
$D^2$ & $2.0321 \times 10^{-9} \, \text{m}^2/\text{s}$ & Diffusion coefficient ($\mathrm{Cl}^-$) \\		
$K_B T$ & $4.14 \times 10^{-21} \, \text{J}$ & Boltzmann energy \\		
$\epsilon_0$ & $8.85 \times 10^{-12} \, \text{C}^2/(\text{N·m}^2)$ & Vacuum permittivity \\		
$e$ & $1.6 \times 10^{-19} \, \text{C}$ & Elementary charge \\		
\hline	
\end{tabular}
\end{table}

\begin{figure}[H]
\begin{subfigure}[b]{0.5\linewidth}%
	\centering
	\includegraphics[height=5cm,    width=5.5cm]{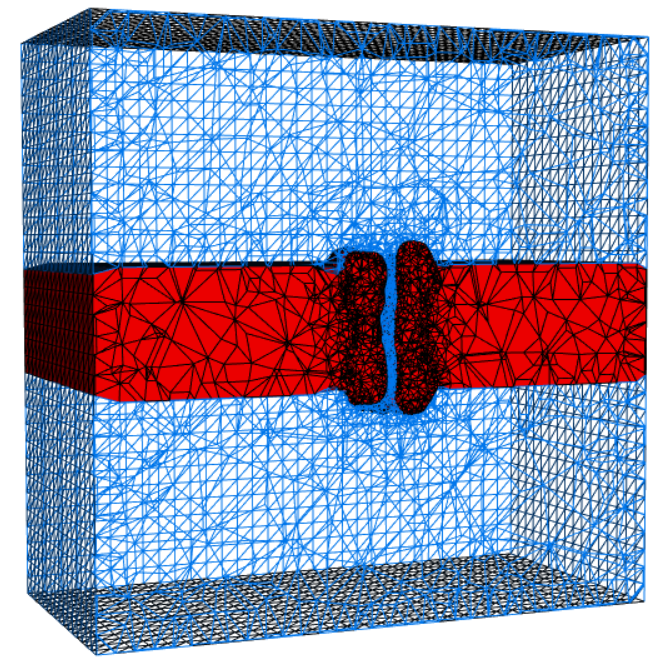}
	\caption{The unstructured tetrahedral volume mesh for the GA channel}
	\label{fig:GA1}
\end{subfigure}
\begin{subfigure}[b]{0.5\linewidth}
	\centering
	\includegraphics[height=5cm,    width=5cm]{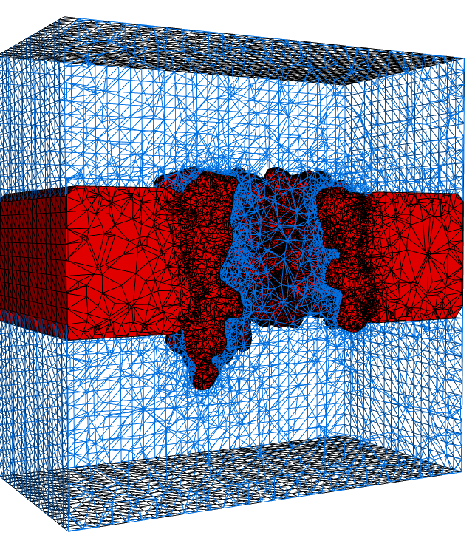}
	\caption{The unstructured tetrahedral volume mesh for the Fhua artificial channel}
	\label{fig:KcsA1}
\end{subfigure}
\caption{\color{black}Training channel.}
\label{fig:trainingchannel}
\end{figure}

The following introduces the two types of channels used for training in Algorithm \ref{RFR}, whose structures are shown in  Fig. \ref {fig:trainingchannel}. On the left is the native ion channel GA, a classic helical dimer protein that has become the most commonly used ion channel model due to its simple structure and well-characterized function. On the right is the engineered ion channel FhuA, which is constructed based on the E. coli FhuA protein (PDB ID: 1BY3) using Algorithm \ref{artificial channel generation}. FhuA was selected as the protein scaffold mainly for the following advantages.
First, as a monomeric protein, it avoids the complexity of subunit assembly encountered in modifying multimeric proteins such as $\alpha$-HL and OmpF. Second, its stable 
$\beta$-barrel structure provides excellent engineering flexibility. In addition, its large native pore diameter offers sufficient potential for further structural modification \cite{Mohammad2010}.

\textcolor{black}{Based on the two training channels(GA and the Fhua artificial channel), we further generated a dataset for model training. The dataset consists of a total of 500 samples,} covering a range of parameter combinations, including voltage, membrane charge density, initial ion concentration, and Gummel iteration relaxation parameters. The specific value ranges and distributions of each parameter are shown in Table~\ref{tab:input_parameters}.

\begin{table}[h]
\centering
\caption{Training set for RFR model}
\begin{tabular}{lcccc}
\hline
\textbf{Input Parameters} & \textbf{Min} & \textbf{Average} & \textbf{Max} & \textbf{mark}\\
\hline
Voltage & 0 & 3.64 & 10 & $\Phi$ in (\ref{2.2})\\

Membrane Charge Density & 0 & 0.28 & 0.9 & $\rho_f$ in (\ref{2.1})\\

Initial Ion Concentration & 0 & 0.29 & 0.9 & $p$ and $n$ in (\ref{2.2})\\

Relaxation Parameter & 0.04 & 0.51 & 0.91 & $\alpha$ in Algorithm \ref{Gummel}\\

Number of Iterations & 12 & 91.5 & 1000 & \makecell[c]{The number of iterations \\obtained from Algorithm \ref{Gummel}}\\

\makecell[c]{Geometric feature vector \\$G$ of channel} & \multicolumn{3}{c}{Determined by Algorithm \ref{for Simple Channel} / \ref{for Complex Channel} } \\
\hline
\end{tabular}

\label{tab:input_parameters}
\end{table}

\begin{figure}[H]
	\centering
	\begin{subfigure}[b]{0.3\linewidth}
		\centering
		\includegraphics[height=4.9cm,    width=2.7cm]{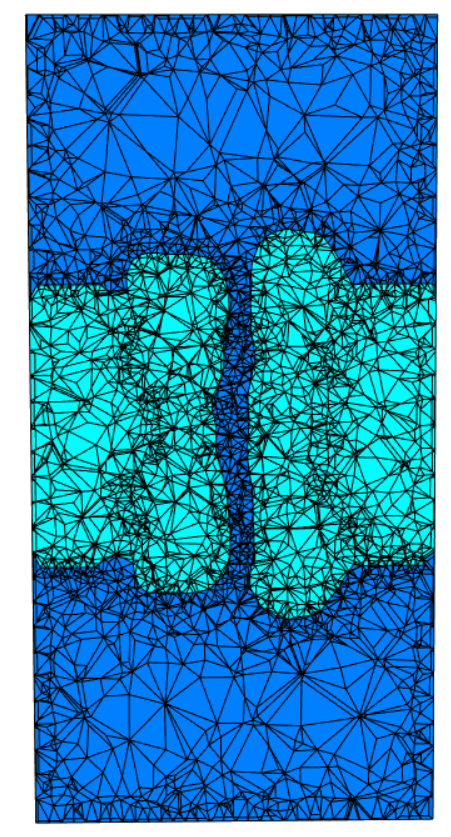} 
		\caption{GA channel cross-section}
		\label{fig:GAcross-section}
	\end{subfigure}
	\hspace{0cm}  
	\begin{subfigure}[b]{0.3\linewidth}
		\centering
		\includegraphics[height=5cm,    width=5cm]{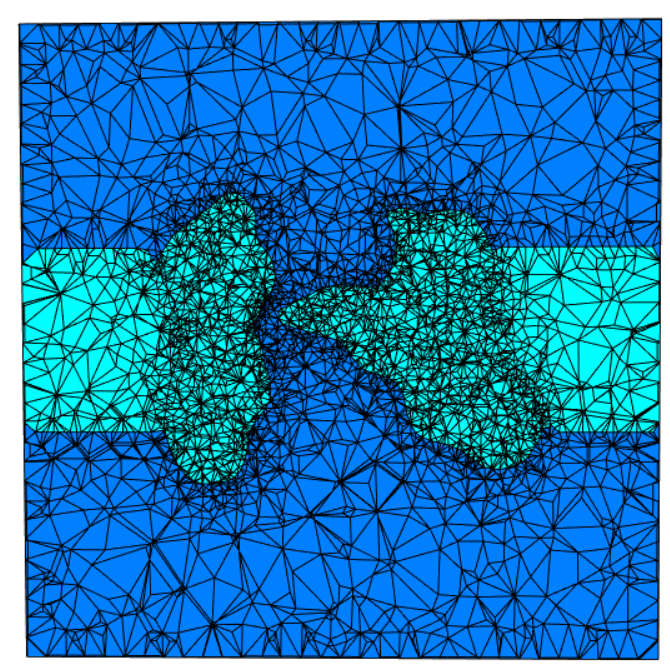} 
		\caption{VDAC channel cross-section}
		\label{fig:VDACcross-section}
		\end{subfigure}
	\hspace{0.5cm}  
	\begin{subfigure}[b]{0.3\linewidth}
		\centering
		\includegraphics[height=5cm,    width=5cm]{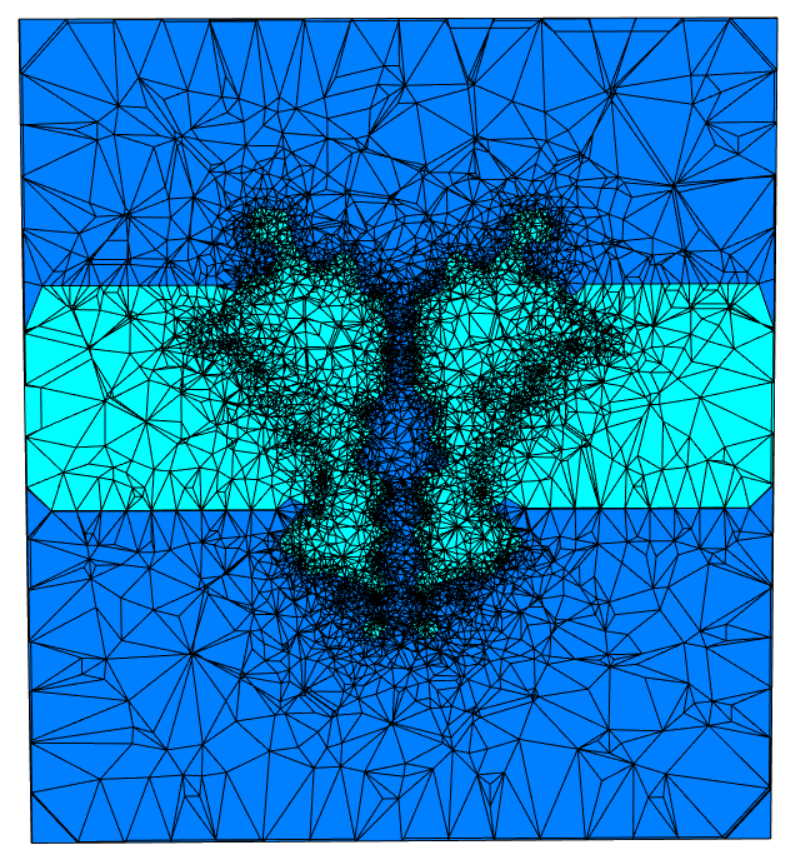} 
		\caption{KcsA channel cross-section}
		\label{fig:KcsAcross-section}
	\end{subfigure}
	\caption{test channel cross-section}
	\label{fig:maincross-sectionfigure}
\end{figure}

The simulations performed on three test channels: GA[Protein Data Bank (PDB) ID: 1MAG], VDAC[PDB ID: 2JK4] and KcsA[PDB ID: 1BL8]. Test Channel 1 is the GA channel, whose cross-sectional mesh is illustrated in Fig. \ref{fig:GAcross-section}. Notably, although the GA channel is also used during model training, the dataset employed for testing is entirely distinct from the training dataset. Test Channel 2 is the VDAC, which features a characteristic $\beta$-barrel structure \cite{Mannella1992}, with its cross-sectional mesh shown in Fig. \ref{fig:VDACcross-section}. The final test channel is the KcsA channel, which exhibits a more complex spatial conformation, including a symmetric arrangement of four subunits, a narrow selectivity filter, and multiple ligand-binding sites. Its cross-sectional view is presented in Fig. \ref{fig:KcsAcross-section}.
By comparing Figs. \ref{fig:maincross-sectionfigure} and \ref{fig:trainingchannel}, it can be seen that the geometries of the VDAC and KcsA channels differ substantially from those of the training channels, presenting a considerable challenge for validating the effectiveness of Algorithm \ref{RFR}.


\subsection{Test channel 1: GA}

In this subsection, Algorithm \ref{RFR} is adopted to solve the PNP system \eqref{2.1}-\eqref{2.2} for simulating the GA channel. It mainly consists of three parts. First, the discrete mesh of the computational domain is introduced. Second, the accuracy of Algorithm \ref{RFR} is verified by comparing its numerical current with experimental and existing simulation data. Finally, the performance of Algorithm \ref{RFR} and Algorithm \ref{Gummel} is compared from the perspective of iteration efficiency.

\begin{figure}[ht]
	\begin{subfigure}[t]{0.5\linewidth}
		\centering
		\includegraphics[height=4.5cm,width=4cm]{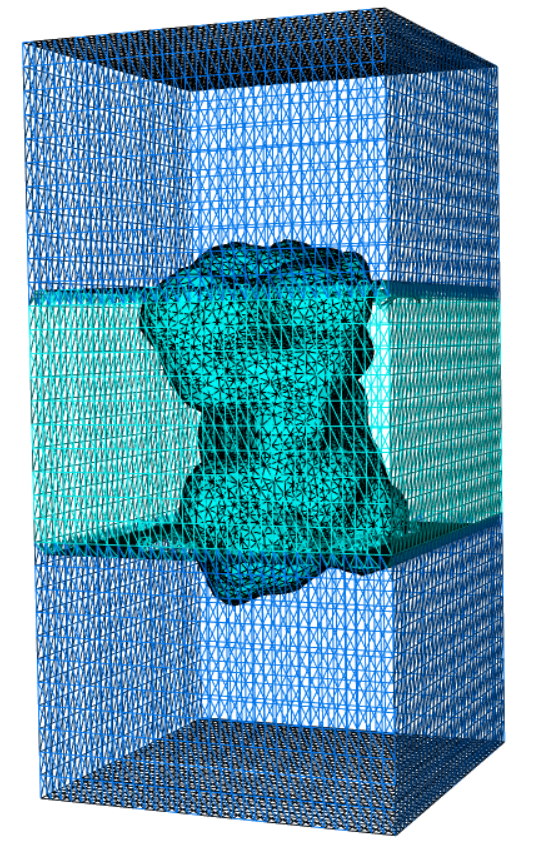}
		\caption{Computational region}  
		\label{fig:GAa}
	\end{subfigure}
	\begin{subfigure}[t]{0.5\linewidth}
		\centering
		\includegraphics[height=3.5cm,width=4.8cm]{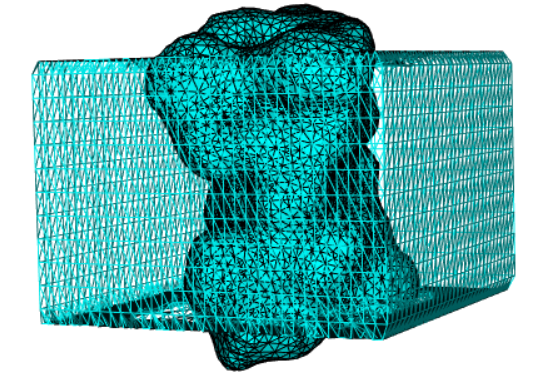}
		\caption{Membrane and protein region}  
		\label{fig:GAb}
	\end{subfigure}
	\begin{subfigure}[t]{0.5\linewidth}
		\centering
		\includegraphics[height=4.5cm,width=4cm]{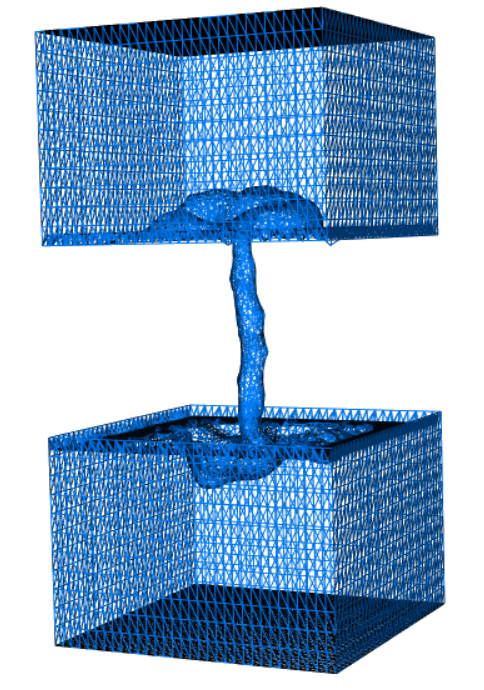}
		\caption{Bulk region}  
		\label{fig:GAc}
	\end{subfigure}
	\begin{subfigure}[t]{0.5\linewidth}
		\centering
		\includegraphics[height=4.5cm,width=3.7cm]{GAchannelcross-section.png}
		\caption{Tetrahedral mesh cross-section}  
		\label{fig:GAd}
	\end{subfigure}
	
	\caption{The unstructured mesh for the GA channel}
	\label{fig:GA}
\end{figure}
\begin{figure}[!htbp]	
\begin{subfigure}[htb]{0.5\linewidth}		
\centering		
\includegraphics[height=5cm,width=4cm]{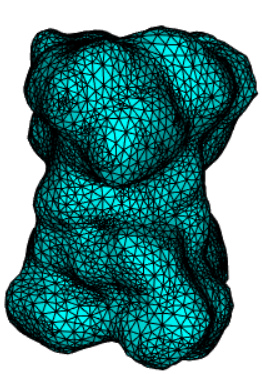}		
\caption{Top view}  
\label{fig:GAiona}	
\end{subfigure}	
\begin{subfigure}[htb]{0.5\linewidth}		
\centering		
\includegraphics[height=4.5cm,width=4.8cm]{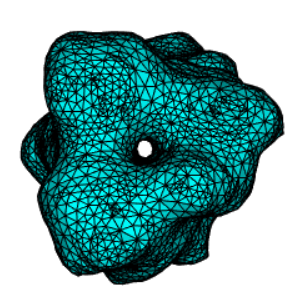}		
\caption{Lateral view}  
 \label{fig:GAionb}	
 \end{subfigure}

 \caption{Triangular mesh of the GA channel surface}	
 \label{fig:GAion}
\end{figure}

Consider the PNP system \eqref{2.1}-\eqref{2.2}, where the solution domain $\Omega$ of the GA channel is defined as $\Omega = [x, y, z] = [-15, 15] \times [-15, 15] \times [-30, 30]$. The discrete mesh is generated using TMSmesh software, with a total of 29092 nodes and 163325 elements. Fig. \ref{fig:GA} illustrates the meshes of the overall and partial regions of the GA domain, revealing that the interface between the membrane and the bulk region is highly irregular. Fig. \ref{fig:GAion} presents the top view and lateral view meshes of the ion channel region, which shows that the channel surface also exhibits significant irregularities.


Subsequently, we employ Algorithm \ref{RFR} to solve the PNP problem \eqref{2.1}-\eqref{2.2} for GA channel simulations. As analytical solutions to the PNP problem \eqref{2.1}-\eqref{2.2} are unavailable, we validate the accuracy and effectiveness of Algorithm \ref{RFR} via comparisons between its numerical results and experimental data as well as existing numerical results. Following common practices for GA channel numerical simulations, we conduct comparisons in terms of ionic current. From the computed PNP solutions, the channel ionic current is further calculated and compared with available experimental and simulation results.


Table \ref{current} presents a comparison of ionic currents among experimental data in \cite{COLE20021974}, numerical results from Algorithm \ref{RFR}, Algorithm \ref{Gummel} and \cite{Tu2013Aparallelfinite}. The numerical current is calculated as the net charge flux through any channel cross-section $\Omega_I$\cite{Tu2013Aparallelfinite, luzh11}:
\begin{equation} \label{curr1}
    I = e \int_{\Omega_I} (J_1 - J_2) dS,
\end{equation}
where the ionic fluxes $J_i$ are given by $$J_i = -\nabla \cdot D_i \left( \nabla p^i + \frac{e}{K_B T} q^i p^i \nabla \phi \right), ~i=1,2$$ 
with $(p^i,\phi)$ denoting the numerical solution to the PNP problem \eqref{2.1}-\eqref{2.2}. 

\textcolor{black}{The results in Table \ref{current} show that ionic current results from Algorithm \ref{RFR} agree well with experimental data in \cite{COLE20021974}, with a smaller standard deviation than the existing simulation result in \cite{Tu2013Aparallelfinite}, implying lower simulation error. This demonstrates that Algorithm \ref{RFR} is accurate and effective for solving the PNP equations in GA channel simulations. Furthermore, the ionic current obtained by the proposed Algorithm \ref{RFR} agrees well with that of the classic Gummel iteration, i.e., Algorithm \ref{Gummel}, indicating that the artificial channel generation method in Algorithm \ref{artificial channel generation} and the geometric feature extraction processes in Algorithms \ref{for Simple Channel}–\ref{for Complex Channel} do not substantially affect the numerical solution.}



\begin{table}[ht]
\centering
\caption{{Ionic current comparison and standard deviation (SD) vs. experimental data}}
\begin{tabular}{cccccccc}   
\toprule   
\makecell[c]{Concentration\\(M)} & \makecell[c]{Voltage\\(mV)} & \makecell[c]{Experimental\\current(pA)\cite{COLE20021974}} & \makecell[c]{current(pA)\\ of \cite{Tu2013Aparallelfinite}} & \makecell[c]{Alg. \ref{RFR}\\current\\(pA)} &  \makecell[c]{Alg. \ref{Gummel}\\current\\(pA)} &\makecell[c]{SD of \\Alg. \ref{RFR}} & \makecell[c]{SD\\ in \cite{Tu2013Aparallelfinite}}\\   
\midrule   
0.1 & 50 & 0.65 & 0.52 & 0.63 & 0.63 & 0.02 & 0.12\\ 
    & 100 & 1.20 & 1.03 & 1.21 & 1.21 & \\ 
    & 150 & 1.71 & 1.56 & 1.73 & 1.73 & \\ 
    & 200 & 2.12 & 2.12 & 2.20 & 2.20 & \\
0.2 & 50 & 1.06 & 0.81 & 0.93 & 0.93 & 0.08 & 0.20\\ 
    & 100 & 1.89 & 1.62 & 1.82 & 1.82 &\\ 
    & 150 & 2.72 & 2.48 & 2.66 & 2.66 & \\ 
    & 200 & 3.51 & 3.39 & 3.47 & 3.47 & \\
0.5 & 50 & 1.66 & 1.31 & 1.57 & 1.57 & 0.21 & 0.69\\ 
    & 100 & 3.46 & 2.64 & 3.15 & 3.15 & \\   
    & 150 & 4.94 & 4.07 & 4.73 & 4.73 &\\ 
    & 200 & 6.55 & 5.61 & 6.32 & 6.32 & \\

\bottomrule
\end{tabular}

\label{current}
\end{table}

\begin{figure}[!ht]
	\begin{subfigure}[t]{0.5\linewidth}
		\centering
\includegraphics[height=5.8cm,width=4.5cm]{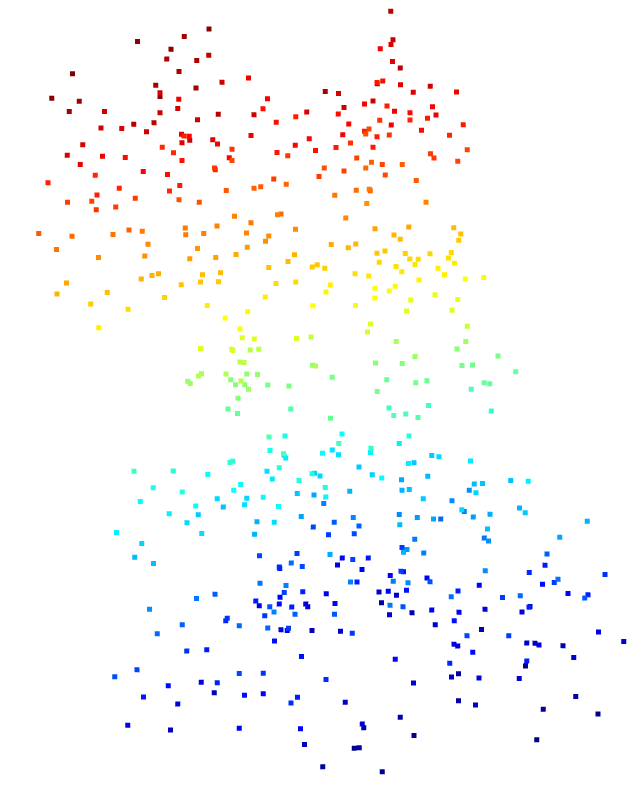}
		\caption{Front view, original GA channel}
		\label{fig:GAatom:a}
	\end{subfigure}
	\begin{subfigure}[t]{0.5\linewidth}
		\centering
\includegraphics[height=5.8cm,width=4cm]{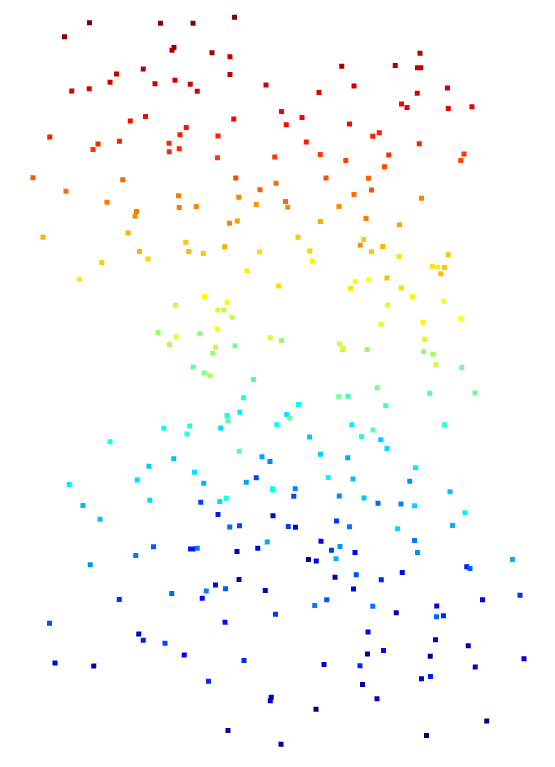}
		\caption{Front view, after Algorithm \ref{for Simple Channel}}
		\label{fig:arg2GA:c}
	\end{subfigure}
	\begin{subfigure}[t]{0.5\linewidth}
		\centering
	\includegraphics[height=4.8cm, width=4.5cm]{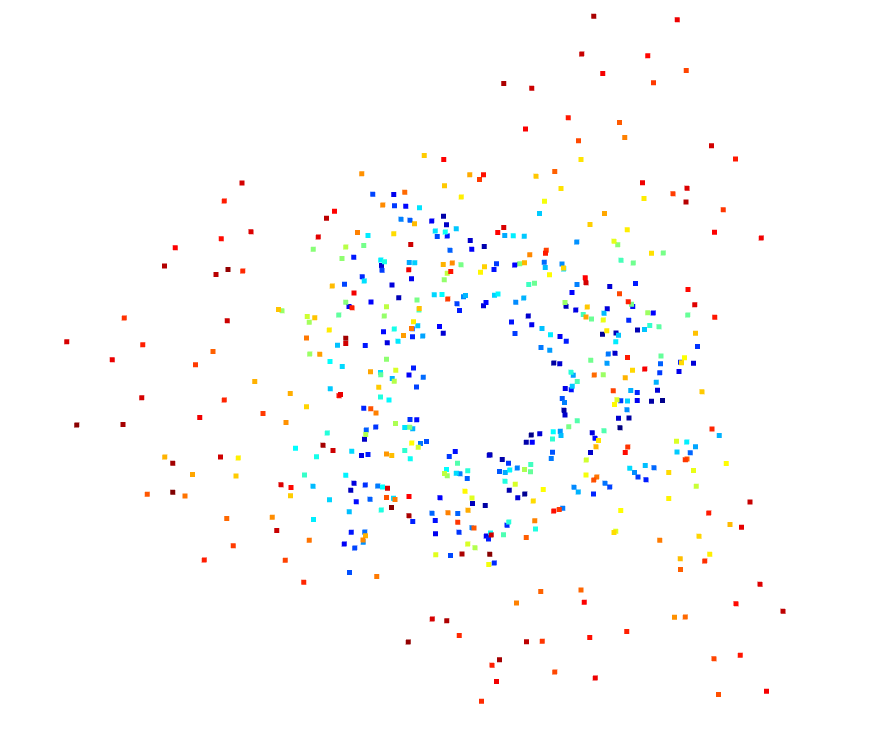}
		\caption{Top view, original GA channel}
 		\label{fig:GAatom:b}
		\label{fig:arg2GA:c}
	\end{subfigure}
	\begin{subfigure}[t]{0.5\linewidth}
		\centering
		\includegraphics[height=4cm,width=4cm]{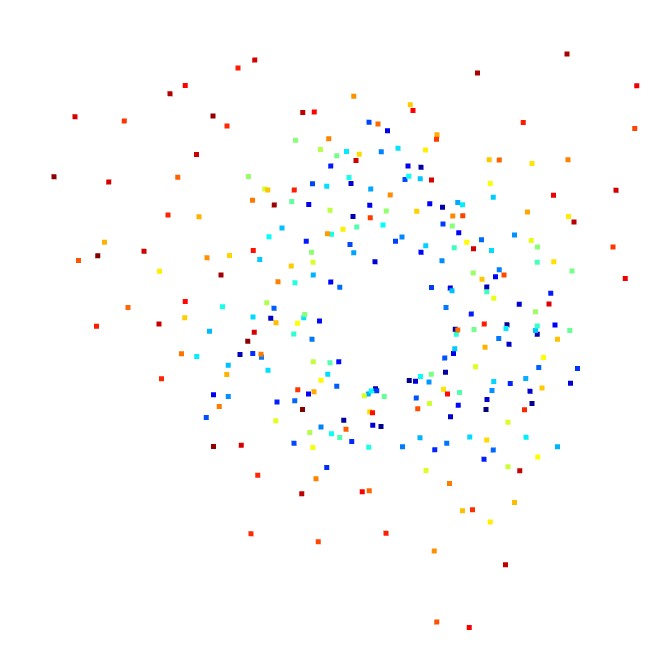}
		\caption{Top view, after Algorithm \ref{for Simple Channel}}
		\label{fig:arg2GA:d}
	\end{subfigure}
	
	\caption{Comparison of the GA Channel before and after the application of the Algorithm \ref{for Simple Channel}}
	\label{fig:arg2GA}
\end{figure}

We next present the performance of Algorithm \ref{RFR} in improving iteration efficiency. Since Algorithm \ref{RFR} adopts Algorithm \ref{for Simple Channel} in the prediction process (see step 3), we first observe the performance of Algorithm \ref{for Simple Channel}. Afterwards, we show the iteration counts of Algorithm \ref{RFR} and compare them with those of the traditional algorithm.

Fig. \ref{fig:arg2GA} shows the point-charge distributions of the GA channel before and after the application of  Algorithm \ref{for Simple Channel}. In the original PQR files of GA, a total of 552 atoms are included. After applying Algorithm \ref{for Simple Channel}, only 200 atoms are retained. By comparing the front and top views of point-charge distributions in Fig. \ref{fig:arg2GA}, it can be observed that the algorithm significantly reduces the number of atoms of GA channel while preserving its geometric characteristics. This shows Algorithm \ref{for Simple Channel} can efficiently capture geometric features of the GA channel, and the reduced atom count can enhances the prediction efficiency of Algorithm \ref{RFR}.

Next, we compare the iteration efficiency between the proposed Algorithm \ref{RFR} and the classical Algorithm \ref{Gummel}. The primary difference between the two algorithms lies in the selection of relaxation parameters. The proposed Algorithm \ref{RFR} adopts relaxation parameters automatically predicted by machine learning, whereas the classical Algorithm \ref{Gummel} uses manually‑selected empirical relaxation parameters. Therefore, to verify the iteration efficiency of Algorithm \ref{RFR}, we select several representative manual relaxation parameters within a reasonable range for Algorithm \ref{Gummel}, and compare its iteration performance with that of Algorithm \ref{RFR}. 

\begin{figure}[ht]
	\begin{subfigure}[t]{0.5\linewidth}
		\centering
		\includegraphics[height=7cm,width=8cm]{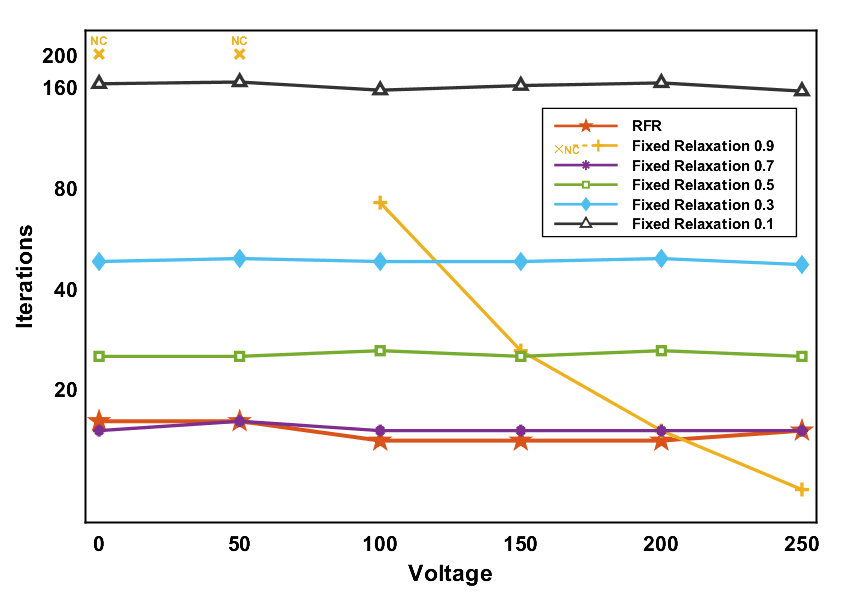}
		\caption{Number of iterations at different voltages}
		\label{fig:GAEiter:a}
	\end{subfigure}
	\begin{subfigure}[t]{0.5\linewidth}
		\centering
		\includegraphics[height=7cm,width=8cm]{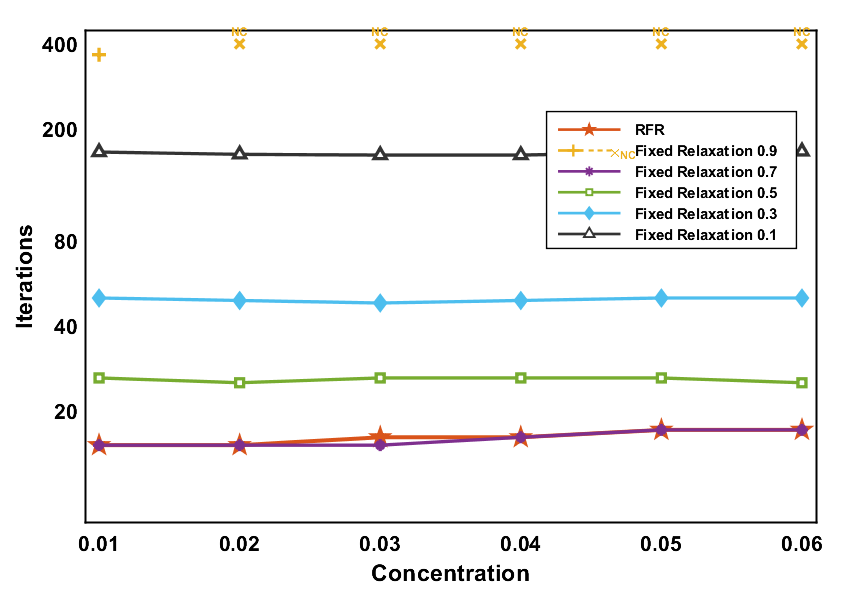}
		\caption{Number of iterations at different concentrations}
		\label{fig:GAninfiter:b}
	\end{subfigure}
	\caption{The number of Gummel iterations}
	\label{fig:GAiter}
\end{figure}

Fig \ref{fig:GAiter} shows the number of Gummel iterations when different relaxation parameters are selected under the circumstances of different voltages and concentrations. The red line represents the numbers of iterations for Algorithm \ref{RFR}, where the relaxation parameters selected by RFR, while the other colored lines indicate the numbers of iterations of Algorithm \ref{Gummel} with manually selected relaxation parameters. As observed from Figure \ref{fig:GAEiter:a}, Algorithm \ref{RFR} based on machine learning achieves nearly the fewest iteration steps under various voltages and concentrations. This demonstrates that the relaxation parameters predicted by Algorithm \ref{RFR} are quasi‑optimal.


Table \ref{Tab:VDACb} presents the number of iteration and CPU times of Algorithm \ref{RFR} and Algorithm \ref{Gummel} for selecting relaxation parameters under various conditions. Algorithm \ref{RFR} employs the RFR method for automatic selection, while Algorithm \ref{Gummel} adopts a manual selection approach. Numerical results show that the number of iterations and CPU time of Algorithm \ref{RFR} are roughly comparable to those of Algorithm \ref{Gummel} under the optimal manually-selected relaxation parameter. Although Algorithm \ref{RFR} incorporates the training and prediction steps of the RFR method compared with Algorithm \ref{Gummel}, its overall CPU time does not increase significantly, indicating that the computational cost for training and prediction of RFR is low. It is worth noting that Algorithm \ref{RFR} exhibits favorable stability. Its number of iterations remains close to the optimum under different voltages and concentrations. By contrast, Algorithm \ref{Gummel} relying on manually-selected parameters exhibits instability. For instance, when the relaxation parameter is manually set to 0.9, Algorithm \ref{Gummel} achieves the best performance at high voltage (200 mV), yet diverges for all concentration cases.

\begin{table}[h]
	\centering
	\caption{\textbf{Relaxation Parameters, numbers of iterations, and CPU Times of Algorithm \ref{RFR} on the GA Channel Under Various Conditions.}}
    \label{Tab:VDACb}
	\begin{minipage}[t]{0.45\linewidth}
		\centering
		\subcaption{\textbf{Result under different voltages}}
		\label{tab:KcsAE1iterations}
		\begin{tabular}{ccccc}
			\toprule
Volt. & Alg. & \makecell[c]{Relaxation \\ parameter} & Iter. & \makecell[c]{CPU \\ time} \\
\midrule
0 & Alg. \ref{Gummel} & 0.10 & 163 & 250.89s \\
     & Alg. \ref{Gummel} & 0.50 & 25 & 38.65s \\
     & Alg. \ref{Gummel} & 0.70 & 15 & 23.82s \\
    
     & \textbf{Alg. \ref{RFR}} & \textbf{0.67} & \textbf{16} & \textbf{26.88s} \\
\midrule
150 & Alg. \ref{Gummel} & 0.10 & 161 & 285.54s \\
     & Alg. \ref{Gummel} & 0.50 & 25 & 38.76s \\
     & Alg. \ref{Gummel} & 0.70 & 15 & 23.69s \\
     
     & \textbf{Alg. \ref{RFR}} & \textbf{0.72} & \textbf{14} & \textbf{23.95s}\\
\midrule
250 & Alg. \ref{Gummel} & 0.10 & 155 & 252.72s \\
     & Alg. \ref{Gummel} & 0.50 & 25 & 38.84s \\
     & Alg. \ref{Gummel} & 0.70 & 15 & 23.82s \\
     
     & \textbf{Alg. \ref{RFR}} & \textbf{0.71} & \textbf{15} & \textbf{25.47s} \\
			\bottomrule
		\end{tabular}
	\end{minipage}
	\hfill
	\begin{minipage}[t]{0.45\linewidth}
		\centering
		\subcaption{\textbf{Result under different concentrations}}
		\label{tab:KcsANINFiterations}
		\begin{tabular}{ccccc}
			\toprule
Conc. & Alg. & \makecell[c]{Relaxation \\ parameter} & Iter. & \makecell[c]{CPU \\ time} \\
\midrule
0.01 & Alg. \ref{Gummel} & 0.10 & 165 & 270.44s \\
     & Alg. \ref{Gummel} & 0.50 & 26 & 41.72s \\
     & Alg. \ref{Gummel} & 0.70 & 15 & 22.81s \\
     
     & \textbf{Alg. \ref{RFR}} & \textbf{0.70} & \textbf{15} & \textbf{25.50s} \\
\midrule
0.03 & Alg. \ref{Gummel} & 0.10 & 161 & 268.66s \\
     & Alg. \ref{Gummel} & 0.50 & 26 & 41.19s \\
     & Alg. \ref{Gummel} & 0.70 & 15 & 22.96s \\
     
     & \textbf{Alg. \ref{RFR}} & \textbf{0.68} & \textbf{16} & \textbf{27.08s}\\
\midrule
0.06 & Alg. \ref{Gummel} & 0.10 & 166 & 264.22s \\
     & Alg. \ref{Gummel} & 0.50 & 25 & 37.58s \\
     & Alg. \ref{Gummel} & 0.70 &  17 & 26.13s \\
     
     & \textbf{Alg. \ref{RFR}} & \textbf{0.67} & \textbf{17} & \textbf{28.86s} \\
			\bottomrule
		\end{tabular}
	\end{minipage}
\end{table}

The results presented in  Fig. \ref{fig:GAiter} and Table \ref{Tab:VDACb} demonstrate that, as the voltage and concentration increase, the iteration counts of Algorithm \ref{RFR} remain stable without significant fluctuations, whereas Algorithm \ref{Gummel} fails to converge under low voltage conditions. Furthermore, since each iteration of the PNP equations for ion channels incurs substantial computational cost, Algorithm \ref{RFR} still exhibits a clear advantage in total CPU time over Algorithm \ref{Gummel}, despite requiring prediction to obtain the relaxation parameters—the prediction time is extremely short and the number of iterations remains near-optimal.

\subsection{Test Case 2}

In this subsection, we apply Algorithm \ref{RFR} to solve the PNP system \eqref{2.1}-\eqref{2.2} for simulating the second channel--VDAC, which is a typical $\beta$-barrel structure and more complex than GA in structure. Similar to the previous subsection, we first introduce the discrete mesh of the computational domain. Then, the accuracy of Algorithm \ref{RFR} is verified by comparing its numerical current with experimental data. Finally, the algorithm is compared with Algorithm \ref{Gummel} to demonstrate its performance in terms of iteration efficiency.

\begin{figure}[H]
	\begin{subfigure}[t]{0.5\linewidth}
		\centering
		\includegraphics[height=5.5cm, width=6cm]{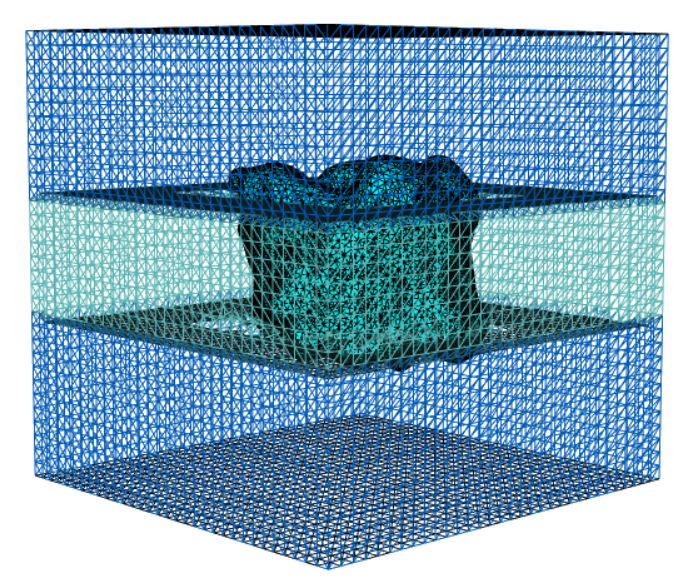}
		\caption{Computational region}
		\label{fig:VDAC:a}
	\end{subfigure}
	\begin{subfigure}[t]{0.5\linewidth}
		\centering
		\includegraphics[height=3.5cm,    width=6cm]{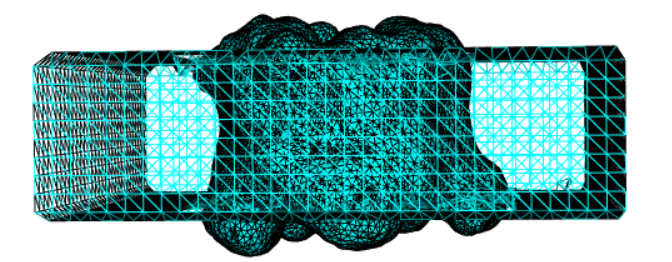}
		\caption{ Membrane and protein region}
		\label{fig:VDAC:b}
	\end{subfigure}
	\begin{subfigure}[t]{0.5\linewidth}
		\centering
		\includegraphics[height=5.5cm, width=5.5cm]{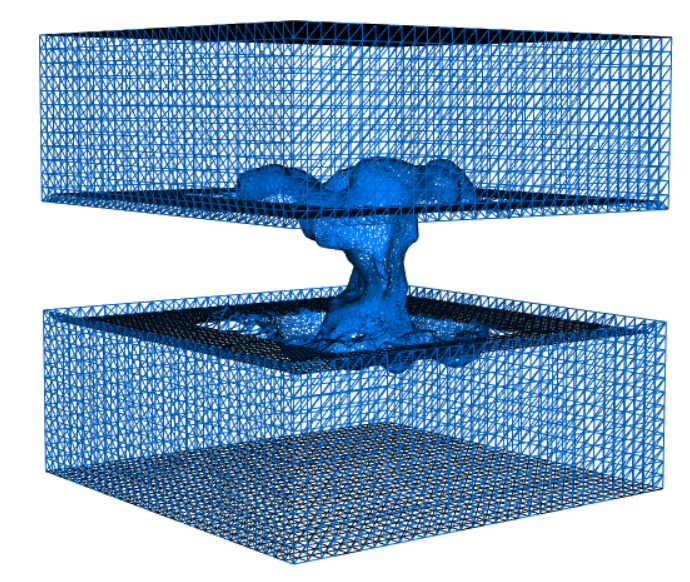}
		\caption{Bulk region}
		\label{fig:VDAC:c}
	\end{subfigure}
	\begin{subfigure}[t]{0.5\linewidth}
		\centering
		\includegraphics[height=5.5cm, width=5.5cm]{VDACchannelcross-section.png}
		\caption{Tetrahedral mesh cross-section}
		\label{fig:VDAC:d}
	\end{subfigure}
	\caption{The unstructured mesh for the VDAC}  
	\label{fig:VDAC}
\end{figure}

Fig. \ref{fig:VDAC:a} presents the entire computational region $\Omega = [x, y, z] = [-20, 70] \times [-45, 45] \times [-40, 50]$, containing total of 93,938 nodes and 569,893 elements. The membrane and protein region $\Omega_{m}$ is marked in green (as shown in Fig. \ref{fig:VDAC:b}), while the remaining part is the bulk solvent region $\Omega_{s}$ (as shown in Fig. \ref{fig:VDAC:c}). Fig. \ref{fig:VDAC:d} displays a cross-sectional view of the entire tetrahedral volume mesh.

\begin{figure}[!htp]	
\begin{subfigure}[t]{0.5\linewidth}		
\centering		
\includegraphics[height=4.5cm,width=4.5cm]{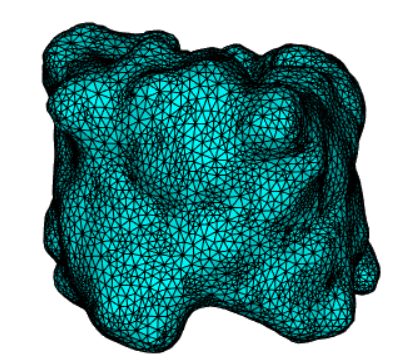}		
\caption{Top view}  
\label{fig:VDACiona}	
\end{subfigure}	
\begin{subfigure}[t]{0.5\linewidth}		
\centering		
\includegraphics[height=4.5cm,width=4.5cm]{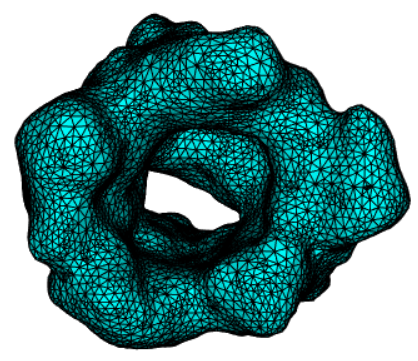}		
\caption{Lateral view}  
\label{fig:VDACionb}	
\end{subfigure}	
\caption{Triangular mesh of the VDAC surface}	
\label{fig:VDAC1}
\end{figure}

In order to better display the shape of VDAC channel, we present the top and lateral view of the surface of VDAC channel, as shown in Fig. \ref{fig:VDAC1}. By comparing Figs. \ref{fig:VDAC1} and \ref{fig:GAion}, it can observed that the surface of the VDAC ion channel is considerably more complex and irregular than that of the GA channel.

\begin{figure}[h] 
\centering
\includegraphics[height=8cm, width=8.5cm]{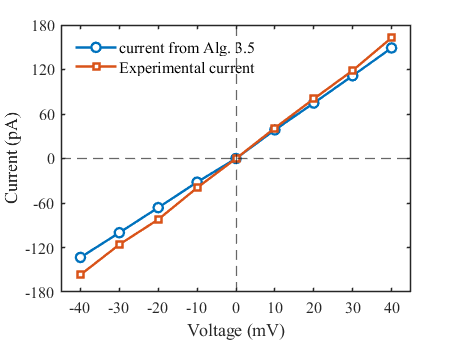} 
\caption{Comparison of the current from Algorithm \ref{RFR} with experimental current \cite{Mertins2012}}
\label{IVcurve3EMN} 
\end{figure}

Due to the lack of corresponding analytical solutions of the PNP equations for VDAC, in order to verify the reliability of Algorithm \ref{RFR}, we verify the correctness of the numerical results by comparing the current 
obtained through Algorithm \ref{RFR} with existing experimental results \cite{Mertins2012}. \textcolor{black}{Fig. \ref{IVcurve3EMN} shows the comparison of currents obtained by Algorithm \ref{RFR} and those from \cite{Mertins2012} for the VDAC channel over the voltage range from -40 mV to 40 mV}, where the numerical currents are are derived from the numerical solutions of the PNP equations, as defined in \eqref{curr1}.  It can be observed that the results of Algorithm \ref{RFR} are in good agreement with the experimental results in \cite{Mertins2012}. 

\begin{figure}[ht]
	\begin{subfigure}[t]{0.5\linewidth}
		\centering
		\includegraphics[height=7cm,    width=5.5cm]{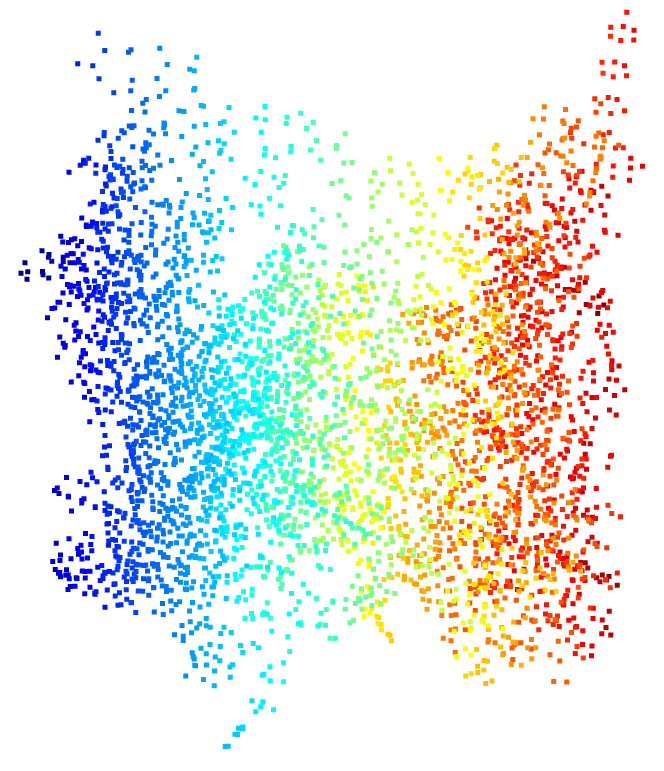}
		\caption{Side view, original VDAC channel}
		\label{fig:VDACatom:a}
	\end{subfigure}
	\begin{subfigure}[t]{0.5\linewidth}
		\centering
\includegraphics[height=7cm,    width=5.5cm]{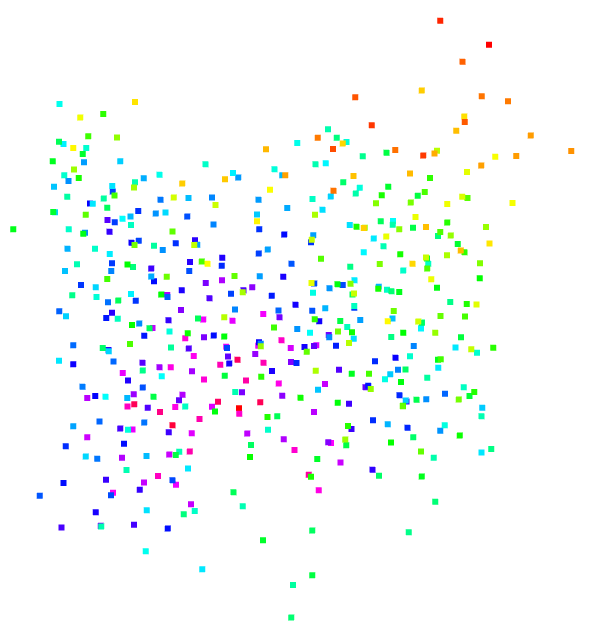}
		\caption{Side view, after Algorithm \ref{for Simple Channel}}
		\label{fig:arg2VDAC:c}
	\end{subfigure}
	\begin{subfigure}[t]{0.5\linewidth}
		\centering
         		\includegraphics[height=5.5cm,    width=5.5cm]{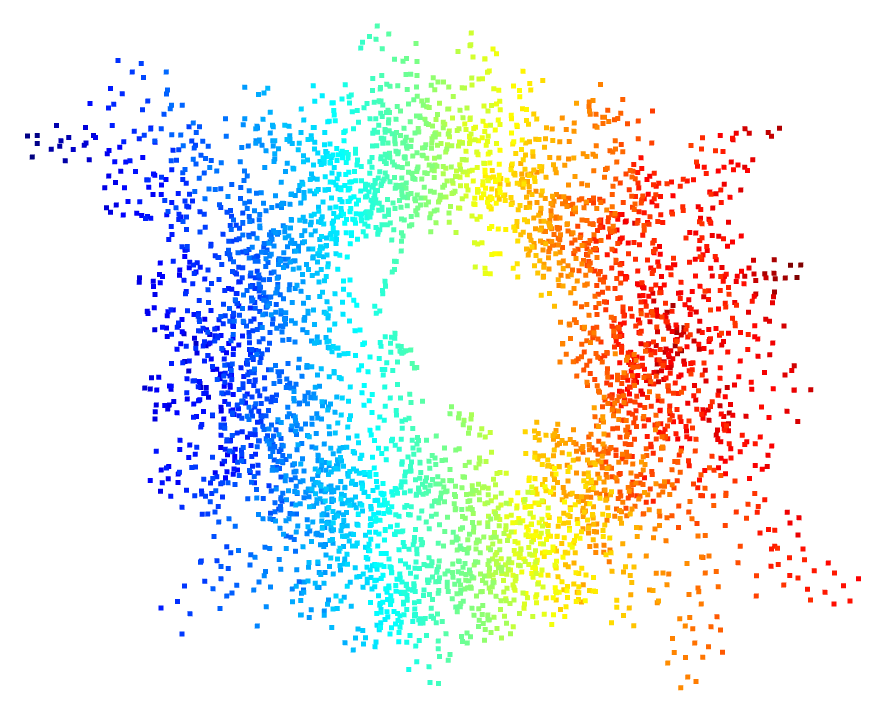}
		\caption{Top view,  original VDAC channel}
		\label{fig:VDACatom:b}
	\end{subfigure}
	\begin{subfigure}[t]{0.5\linewidth}
		\centering
		\includegraphics[height=5.5cm,    width=5.5cm]{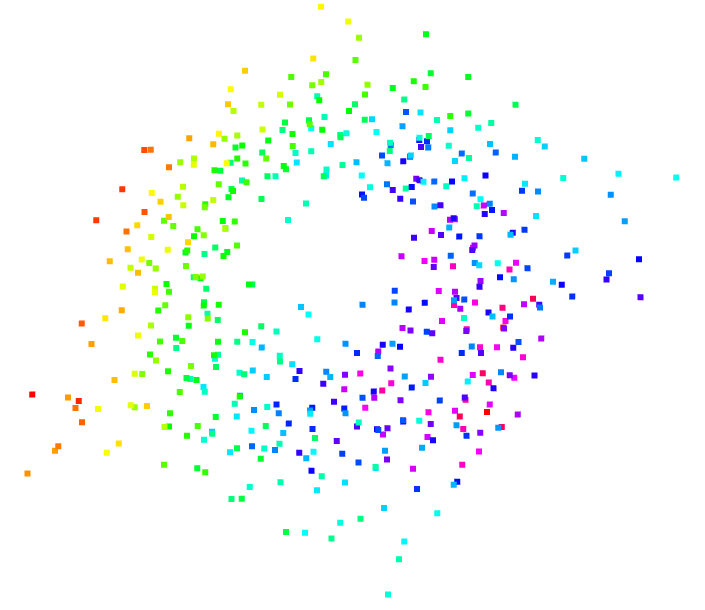}
		\caption{Top view, after Algorithm \ref{for Simple Channel}}
		\label{fig:arg2VDAC:d}
	\end{subfigure}
	\caption{ Comparison of the VDAC Channel before and after the application of Algorithm \ref{for Simple Channel}}
    \label{fig:VDACatom}
\end{figure}

We next demonstrate the performance of Algorithm \ref{RFR} in improving iteration efficiency. Similar to the previous section, we first illustrate the effectiveness of Algorithm \ref {for Simple Channel}, which is utilized in Step 3 of Algorithm \ref{RFR}. Fig. \ref{fig:VDACatom} displays the point-charge distributions of the VDAC channel before and after the application of Algorithm \ref {for Simple Channel}. Panels (a) and (c) present the side view and top view of the original VDAC channel protein with a total of 4394 atoms, respectively. Panels (b) and (d) show the corresponding results after applying Algorithm~\ref{for Simple Channel}, with only 500 atoms remaining. This demonstrates that the algorithm successfully preserves the key geometric features of the channel while significantly reducing the number of point charges. This substantial reduction in the number of point charges enables the RFR model to achieve faster prediction speeds.

\begin{figure}[H]
	\begin{subfigure}[t]{0.5\linewidth}
		\centering
		\includegraphics[height=7cm,    width=8.2cm]{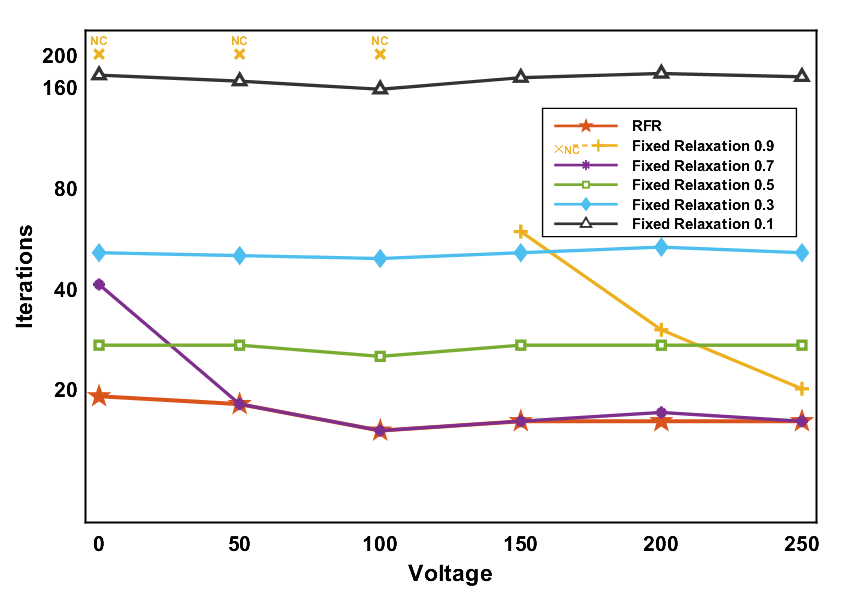}
		\caption{Number of iterations at different voltages}
		\label{fig:VDACEiter:a}
	\end{subfigure}
	\begin{subfigure}[t]{0.5\linewidth}
		\centering
		\includegraphics[height=7cm,    width=8.2cm]{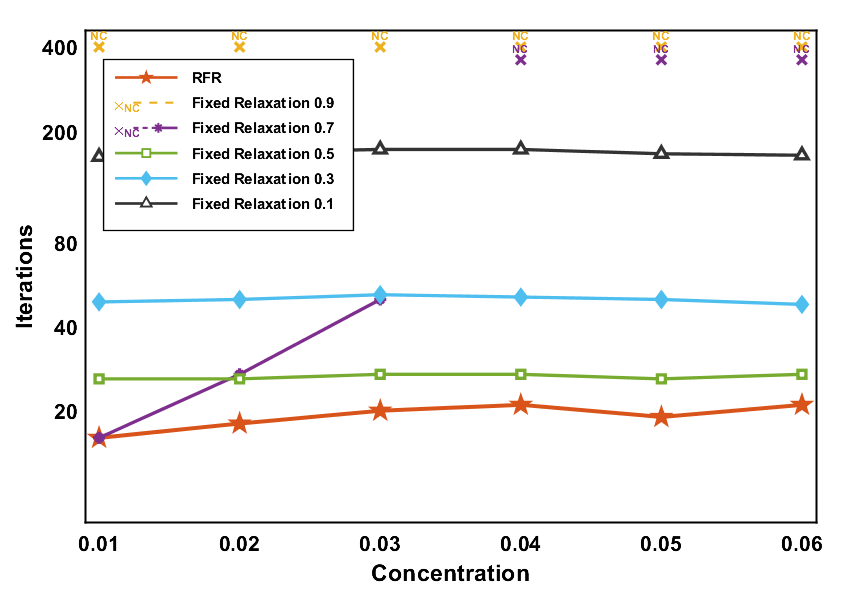}
		\caption{Number of iterations at different concentrations}
		\label{fig:VDACninfiter:b}
	\end{subfigure}
	\caption{Comparison of the number of Gummel iterations between Algorithm \ref{RFR} and Algorithm \ref{Gummel}}
	\label{fig:VDACiter}
\end{figure}

Next, we compare the iteration counts of Algorithms \ref{RFR} and \ref{Gummel}. For Algorithm \ref{RFR}, its relaxation parameter is determined via the RFR method, while the relaxation parameters for Algorithm \ref{Gummel} are selected manually. Fig. \(\ref{fig:VDACiter}\) shows the number of Gummel iterations when different relaxation parameters are selected
under the circumstances of different voltages and concentrations. The red line represents the number of iterations of Algorithm \ref{RFR}, while the other colored lines indicate the number of iterations for Algorithm \ref{Gummel}  with different relaxation parameters. It can be seen that Algorithm \ref{RFR} achieves the smallest number of iterations and exhibits greater stability. In contrast, certain manually selected relaxation parameters lead to instability for different physical parameters.
 For example, when the relaxation parameter is manually set to 0.7 (purple curve), although it achieves near-optimal performance under different voltage conditions, it fails to converge under relatively high concentrations. In contrast, the relaxation parameters determined by the RFR method do not exhibit this issue and maintain a stable number of iterations across all tested conditions. This also demonstrates the effectiveness of Algorithm \ref{RFR}.

Table \ref{Tab:VDACiterb} further compares the CPU time consumption of Algorithm \ref{RFR} and Algorithm \ref{Gummel}. Numerical results show that the number of iterations and CPU time of Algorithm 3.5 are roughly equivalent to those of Algorithm 3.1 under the optimal artificial relaxation parameter. Compared with Algorithm 3.1, although Algorithm 3.5 incorporates the training and prediction stages of the RFR method, its overall CPU time does not increase significantly, indicating low computational overhead for RFR training and prediction. Notably, Algorithm 3.5 exhibits strong stability, with its iteration count consistently close to the optimal value across different voltages and concentrations. In contrast, Algorithm 3.1, which relies on manually selected parameters, suffers from poor stability. For instance, when the manually set relaxation parameter equals 0.7, Algorithm 3.1 achieves optimal performance under various voltage conditions yet diverges at a concentration of 0.06 M.

\begin{table}[h]
	\centering
	\caption{\textbf{Comparison of computational efficiency between Alg. \ref{RFR} and Alg. \ref{Gummel} under different conditions.}}
    \label{Tab:VDACiterb}
	\begin{minipage}[t]{0.45\linewidth}
		\centering
		\subcaption{\textbf{Result under different voltages}}
		\label{tab:KcsAE1iterations}
		\begin{tabular}{ccccc}
			\toprule
Volt. & Alg. & \makecell[c]{Relaxation \\ parameter} & Iter. & \makecell[c]{CPU \\ time} \\
\midrule
0 & Alg. \ref{Gummel} & 0.10 & 173 & 448.77s \\
     & Alg. \ref{Gummel} & 0.50 & 27 & 69.90s \\
     & Alg. \ref{Gummel} & 0.70 & 41  &  109.13s \\
     & \textbf{Alg. \ref{RFR}} & \textbf{0.63} & \textbf{19} & \textbf{53.17s} \\
\midrule
150 & Alg. \ref{Gummel} & 0.10 & 170 & 445.05s \\
     & Alg. \ref{Gummel} & 0.50 & 27 & 70.69s \\
     & Alg. \ref{Gummel} & 0.70 & 16  &  42.76s \\
     & \textbf{Alg. \ref{RFR}} & \textbf{0.71} & \textbf{16} & \textbf{46.25s}\\
\midrule
250 & Alg. \ref{Gummel} & 0.10 & 171 & 436.88s \\
     & Alg. \ref{Gummel} & 0.50 & 27 & 72.38s \\
     & Alg. \ref{Gummel} & 0.70 &  16 & 42.22s \\
     & \textbf{Alg. \ref{RFR}} & \textbf{0.71} & \textbf{16} & \textbf{46.14s} \\
			\bottomrule
		\end{tabular}
	\end{minipage}
	\hfill
	\begin{minipage}[t]{0.45\linewidth}
		\centering
		\subcaption{\textbf{Result under different concentrations}}
		\label{tab:KcsANINFiterations}
		\begin{tabular}{ccccc}
			\toprule
Conc. & Alg. & \makecell[c]{Relaxation \\ parameter} & Iter. &\makecell[c]{CPU \\ time}\\
\midrule
0.01 & Alg. \ref{Gummel} & 0.10 & 162 & 424.34s \\
     & Alg. \ref{Gummel} & 0.50 & 26 & 68.17s \\
     & Alg. \ref{Gummel} & 0.70 & 16 &  43.83s \\
     & \textbf{Alg. \ref{RFR}} & \textbf{0.69} & \textbf{16} & \textbf{47.17s} \\
\midrule
0.03 & Alg. \ref{Gummel} & 0.10 & 172 & 452.85s \\
     & Alg. \ref{Gummel} & 0.50 & 27 & 70.60s \\
     & Alg. \ref{Gummel} & 0.70 &  50 & 137.11s  \\
     & \textbf{Alg. \ref{RFR}} & \textbf{0.61} & \textbf{20} & \textbf{55.31s}\\
\midrule
0.06 & Alg. \ref{Gummel} & 0.10 & 164 & 455.28s \\
     & Alg. \ref{Gummel} & 0.50 & 27 & 71.02s \\
     & Alg. \ref{Gummel} & 0.70 &  - &  - \\
     & \textbf{Alg. \ref{RFR}} & \textbf{0.54} & \textbf{23} & \textbf{63.61s} \\
			\bottomrule
		\end{tabular}
	\end{minipage}
\end{table}

\subsection{Test Case 3}
In this subsection, we apply Algorithm \ref{RFR} to the KcsA channel, which has a more complex spatial conformation, including the symmetric arrangement of four subunits, a narrow selectivity filter, and multiple ligand binding sites. 

\begin{figure}[H]
	\begin{subfigure}[t]{0.5\linewidth}
		\centering
		\includegraphics[height=5.5cm,    width=5.5cm]{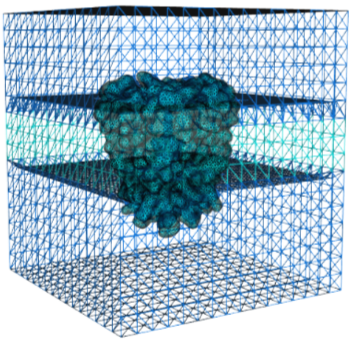}
		\caption{Computational region}
		\label{fig:KcsA:a}
	\end{subfigure}
	\begin{subfigure}[t]{0.5\linewidth}
		\centering
		\includegraphics[height=3cm,    width=5.5cm]{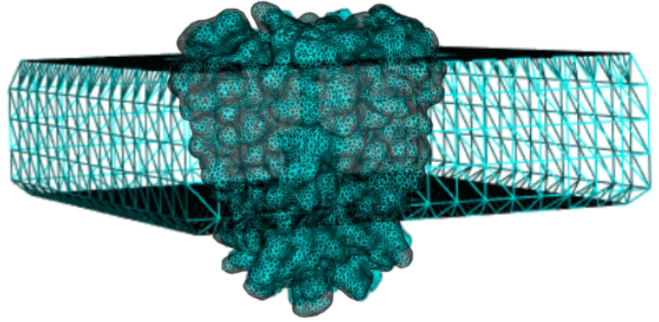}
		\caption{Membrane and protein region}
		\label{fig:KcsA:b}
	\end{subfigure}

	\begin{subfigure}[t]{0.5\linewidth}
		\centering
		\includegraphics[height=5.5cm,    width=5.5cm]{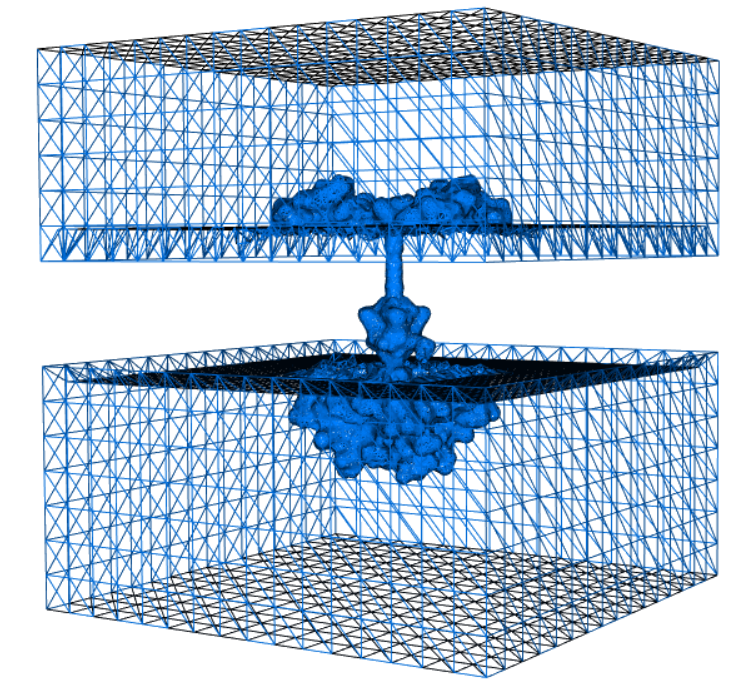}
		\caption{Bulk region}
		\label{fig:KcsA:c}
	\end{subfigure}
	\begin{subfigure}[t]{0.5\linewidth}
		\centering
		\includegraphics[height=5.5cm,    width=5.5cm]{KcsAchannelcross-section.png}
		\caption{Tetrahedral mesh cross-section}
		\label{fig:KcsA:d}
	\end{subfigure}
	\caption{The unstructured mesh for the KcsA channel}
	\label{fig:KcsA}
\end{figure}

The layout of the KcsA channel on the grid is shown in Fig. \ref{fig:KcsA}. Fig. \ref{fig:KcsA1} presents the top and lateral view of the surface of the KcsA channel.  The solution domain $\Omega = [-50, 50] \times [-50, 30] \times [-50, 50] ~\AA$. The mesh of $\Omega$ contains 102,572 nodes and 643,832 elements.

\begin{figure}[H]	
\begin{subfigure}[t]{0.5\linewidth}		
\centering		
\includegraphics[height=5.5cm,width=4.5cm]{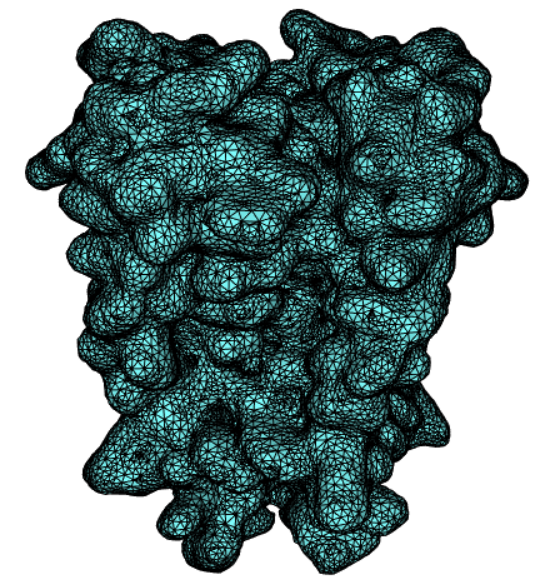}		
\caption{Top view}  
\label{fig:VDACiona}	
\end{subfigure}	
\begin{subfigure}[t]{0.5\linewidth}		
\centering		
\includegraphics[height=4.5cm,width=4.5cm]{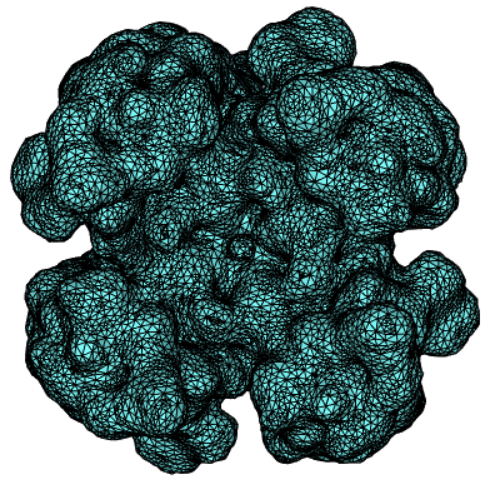}		
\caption{Lateral view}  
\label{fig:VDACionb}	
\end{subfigure}	
\caption{Triangular mesh of the KcsA surface}	
\label{fig:KcsA1}
\end{figure}

\begin{figure}[h] 
\centering
\includegraphics[height=8cm, width=8.5cm]{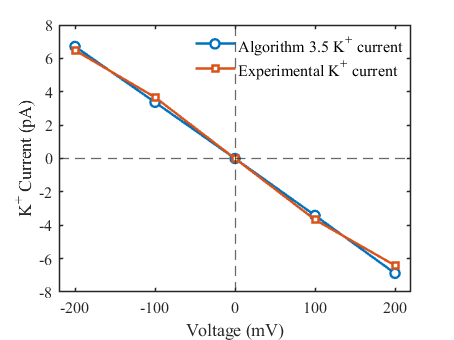} 
\caption{Comparison of the K$^+$ current from Algorithm \ref{RFR} with experimental K$^+$ current}
\label{IVcurve1BL8} 
\end{figure}

Similar to the previous cases, we compare the current obtained by Algorithm \ref{RFR} with experimental results to demonstrate the accuracy of Algorithm \ref{RFR} in solving the PNP equations for the KcsA channel. \textcolor{black}{Fig.~\ref{IVcurve1BL8} shows the K$^+$ ion current of the KcsA channel over the voltage range of $[-200$ mV, $200$ mV]. The red curve denotes the experimental K$^+$ ion current reported in \cite{Liu2017}, while the blue curve corresponds to that calculated using Algorithm \ref{RFR}}. The good agreement between these two sets of curves verifies the computational accuracy of Algorithm \ref{RFR}.

\begin{figure}[H]
	\begin{subfigure}[t]{0.5\linewidth}
		\centering
		\includegraphics[height=5.5cm,    width=5.5cm]{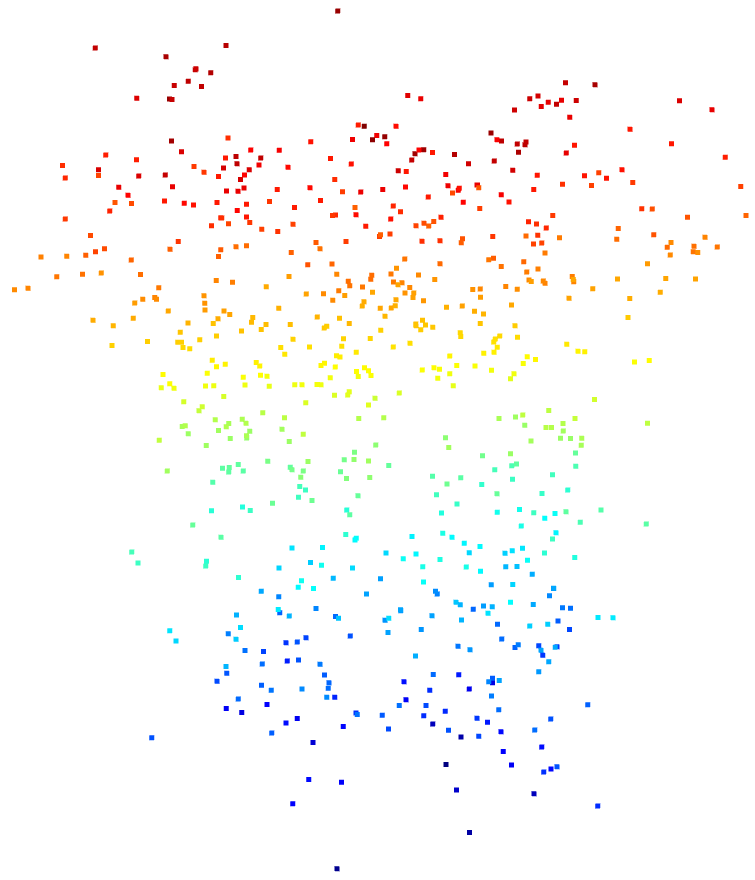}
		\caption{Side view, original KcsA channel}
		\label{fig:KcsAatom:a}
	\end{subfigure}
	\begin{subfigure}[t]{0.5\linewidth}
		\centering
        \includegraphics[height=5.5cm,    width=5.5cm]{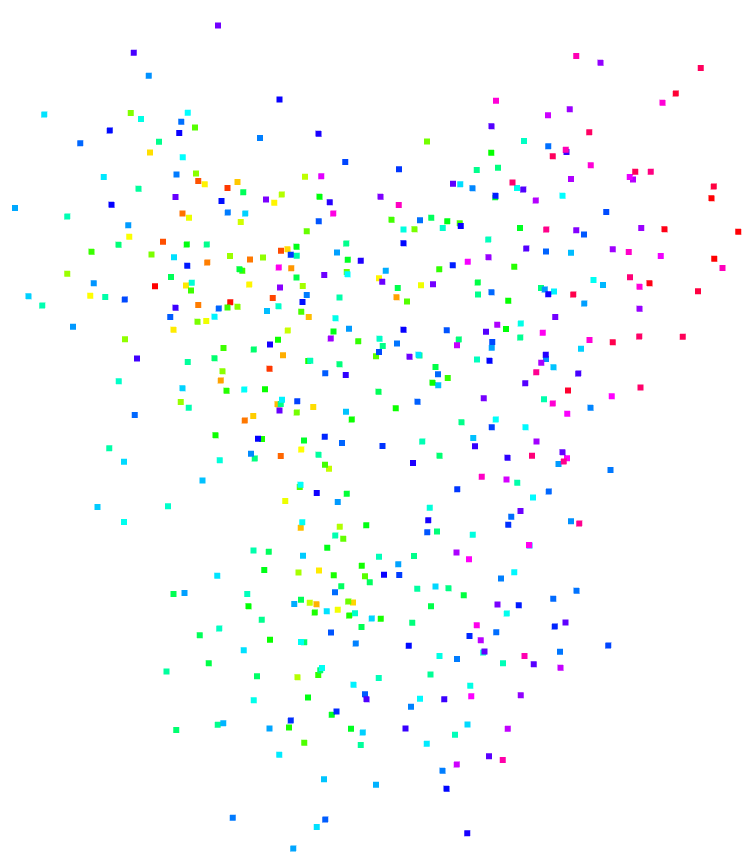}
		\caption{Side view, after Algorithm \ref{for Complex Channel}}
		\label{fig:arg2KcsA:c}
	\end{subfigure}
 
	\begin{subfigure}[t]{0.5\linewidth}
		\centering
        	 \includegraphics[height=5cm,    width=5cm]{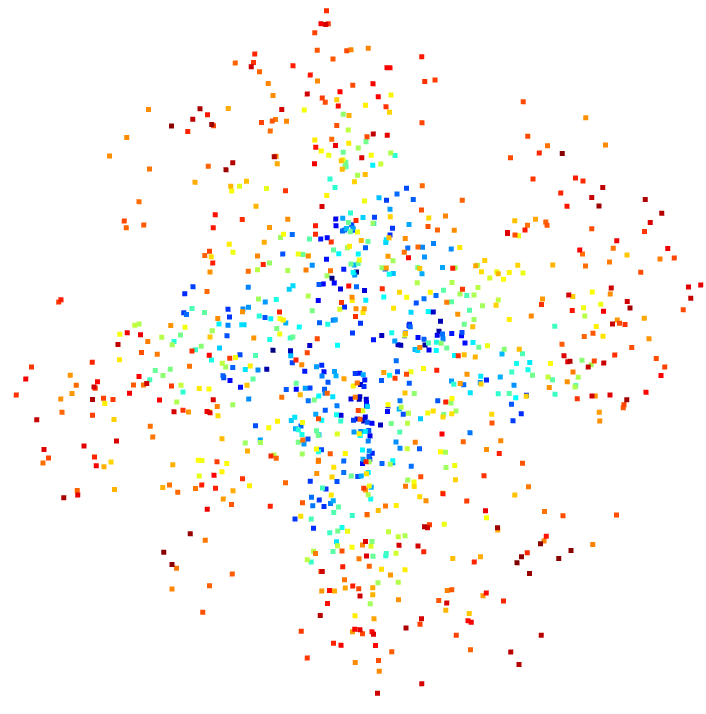}
		 \caption{Top view, original KcsA channel}
		 \label{fig:KcsAatom:b}
	\end{subfigure}
	\begin{subfigure}[t]{0.5\linewidth}
		\centering
		\includegraphics[height=5cm,    width=5cm]{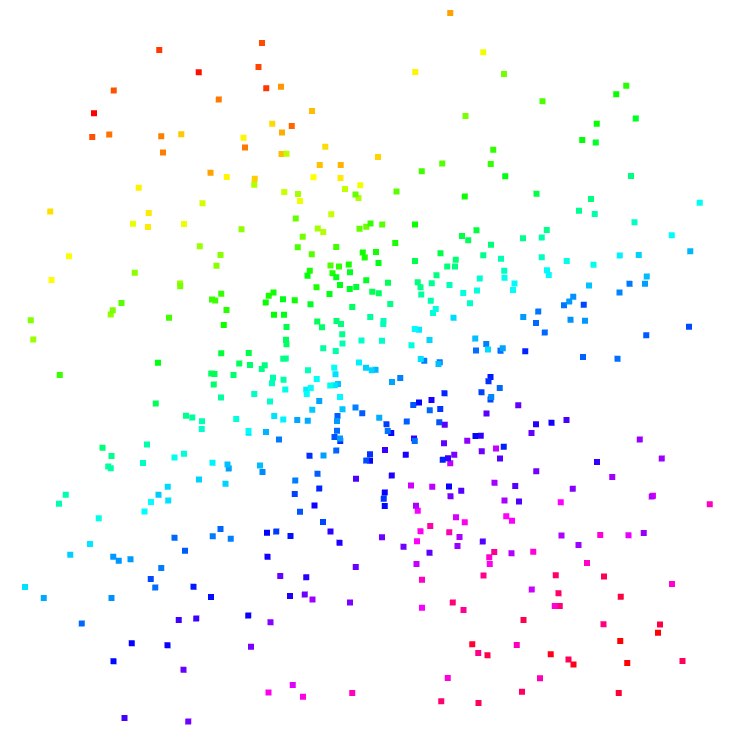}
		\caption{Top view, after Algorithm  \ref{for Complex Channel}}
		\label{fig:arg2KcsA:d}
	\end{subfigure}
	\caption{Comparison of KcsA channel performance before and after Algorithm \ref{for Complex Channel}}
	\label{fig:arg2KcsA}
\end{figure}

The complex geometric features of the KcsA channel  (containing 5995 atoms) make it very difficult to extract its geometric features. Hence, Algorithm \ref{for Simple Channel} can no longer meet the requirements for this extraction task. We adopt Algorithm \ref{for Complex Channel} in Algorithm \ref{RFR} to extract the geometric features of the KcsA channel. Fig. \ref{fig:arg2KcsA} shows the effect of Algorithm \ref{for Complex Channel} with 800 atoms remaining after processing. This demonstrates
that the algorithm preserves the key geometric features of the channel while significantly reducing the number of point charges.

\begin{figure}[H]
	\begin{subfigure}[t]{0.5\linewidth}
		\centering
		\includegraphics[height=7cm,    width=8.2cm]{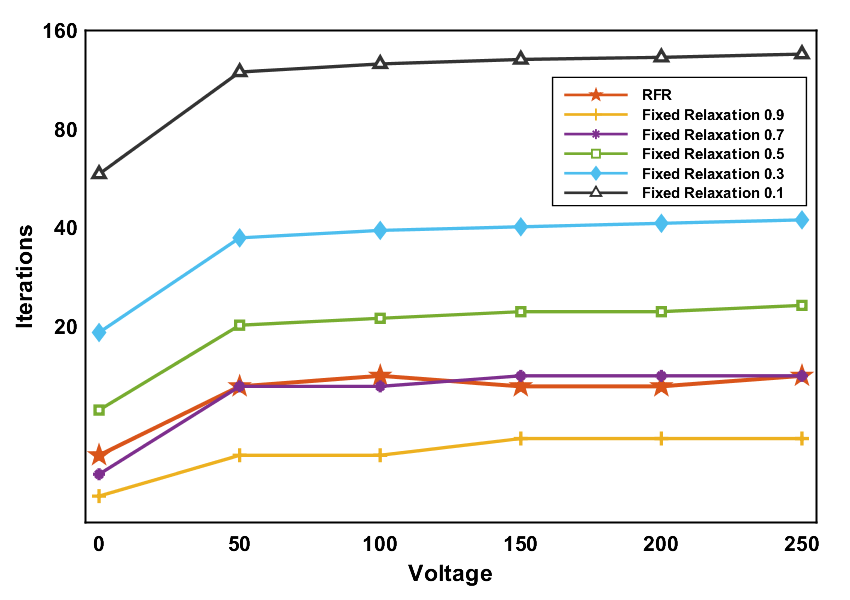}
		\caption{Number of iterations at different voltages}
		\label{fig:KcsAEiter:a}
	\end{subfigure}
	\begin{subfigure}[t]{0.5\linewidth}
		\centering
		\includegraphics[height=7cm,    width=8.2cm]{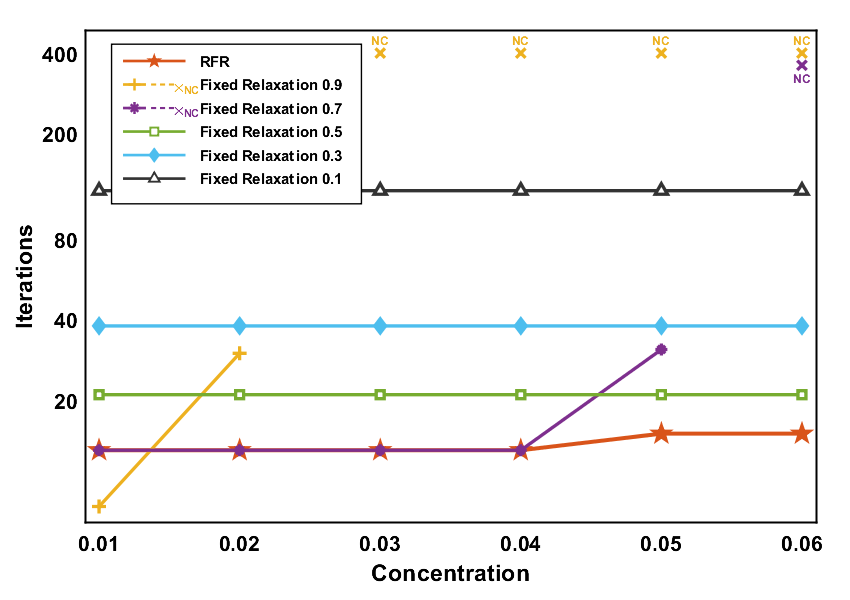}
		\caption{Number of iterations at different concentrations}
		\label{fig:KcsAninfiter:b}
	\end{subfigure}
	\caption{Comparison of the number of Gummel iterations between Algorithm \ref{RFR} and Manual Selection of Relaxation Parameters}
	\label{fig:KcsAiter}
\end{figure}

Fig. \ref{fig:KcsAiter} present the iteration counts of the Gummel iteration for Algorithms \ref{RFR} and \ref{Gummel} on the KcsA channel. The results show that Algorithm \ref{RFR} achieves near-optimality under different applied voltages and is optimal at different given concentrations.  However, the optimal performance of Algorithm \ref{Gummel} relies on manually tuned parameters, which works well for certain conditions but lead to numerical divergence at some high concentrations (e.g., 0.09 M). Therefore, Algorithm 3.5 achievers better robustness and more reliable overall performance.  

Table \ref{Tab:KcsAEiterb} further compares CPU times of Algorithm \ref{RFR} and Algorithm \ref{Gummel}, where manually selected relaxation parameters are 
used for Algorithm \ref{Gummel} on the KcsA channel. 
Algorithm \ref{RFR} exhibits CPU time conparable to Algorithm \ref{Gummel} using optimal parameters. Compared with Algorithm \ref{Gummel}, Algorithm \ref{RFR} adds the training and prediction steps of RFR. This indicates that RFR introduces only minor computational overhead. Together with its near-optimal iteration counts, Algorithm \ref{RFR} achieves stable near-optimal performance with a similar computational cost as Algorithm \ref{Gummel}.

The above figures and tables demonstrate that even on channels with more complex shapes and interfaces, Algorithm \ref {RFR} select the relaxation parameters are still superior. Moreover, under different application conditions, compared with the relaxation parameters manually selected, the number of iterations is significantly more stable. 

\begin{table}[h]
	\centering
	\caption{\textbf{Comparison of computational efficiency between Alg. \ref{RFR} and Alg. \ref{Gummel} under different conditions.}}
    \label{Tab:KcsAEiterb}
	\begin{minipage}[t]{0.45\linewidth}
		\centering
		\subcaption{\textbf{Result under different voltages}}
		\label{tab:KcsAE1iterations}
		\begin{tabular}{ccccc}
			\toprule
Volt. & Alg. & \makecell[c]{Relaxation \\ parameter} & Iter. & \makecell[c]{CPU \\ time} \\
\midrule
0 & Alg. \ref{Gummel} & 0.10 & 58 & 617.26s \\
     & Alg. \ref{Gummel} & 0.50 & 11 & 108.70s \\
     & Alg. \ref{Gummel} & 0.70 & 7 & 66.74s \\
     & \textbf{Alg. \ref{RFR}} & \textbf{0.66} & \textbf{8} & \textbf{80.33s} \\
\midrule
150 & Alg. \ref{Gummel} & 0.10 & 130 & 1406.22s \\
     & Alg. \ref{Gummel} & 0.50 & 22 & 227.09s \\
     & Alg. \ref{Gummel} & 0.70 &  14 & 149.67s \\
     & \textbf{Alg. \ref{RFR}} & \textbf{0.71} & \textbf{13} & \textbf{133.53s}\\
\midrule
250 & Alg. \ref{Gummel} & 0.10 & 135 & 1439.42s \\
     & Alg. \ref{Gummel} & 0.50 & 23 & 236.78s \\
     & Alg. \ref{Gummel} & 0.70 & 14  & 139.34s \\
     & \textbf{Alg. \ref{RFR}} & \textbf{0.71} & \textbf{14} & \textbf{142.22s} \\
			\bottomrule
		\end{tabular}
	\end{minipage}
	\hfill
	\begin{minipage}[t]{0.45\linewidth}
		\centering
		\subcaption{\textbf{Result under different concentrations}}
		\label{tab:KcsANINFiterations}
		\begin{tabular}{ccccc}
			\toprule
Conc. & Alg. & \makecell[c]{Relaxation \\ parameter} & Iter. & \makecell[c]{CPU \\ time} \\
\midrule
0.01 & Alg. \ref{Gummel} & 0.10 & 122 & 1278.67s \\
     & Alg. \ref{Gummel} & 0.50 & 21 & 213.08s \\
     & Alg. \ref{Gummel} & 0.70 &  13 &  128.11s \\
     & \textbf{Alg. \ref{RFR}} & \textbf{0.70} & \textbf{13} & \textbf{134.01s} \\
\midrule
0.03 & Alg. \ref{Gummel} & 0.10 & 122 & 1276.08s \\
     & Alg. \ref{Gummel} & 0.50 & 21 & 212.98s \\
     & Alg. \ref{Gummel} & 0.70 & 13  &  128.02s \\
     & \textbf{Alg. \ref{RFR}} & \textbf{0.71} & \textbf{13} & \textbf{134.99s}\\
\midrule
0.06 & Alg. \ref{Gummel} & 0.10 & 122 & 1280.13s \\
     & Alg. \ref{Gummel} & 0.50 & 21 & 210.32s \\
     & Alg. \ref{Gummel} & 0.70 & -  &  - \\
     & \textbf{Alg. \ref{RFR}} & \textbf{0.63} & \textbf{15} & \textbf{154.25s} \\
			\bottomrule
		\end{tabular}
	\end{minipage}
\end{table}

\section{Conclusion}
This paper presents a method for selecting the relaxation parameters in the Gummel iteration for solving the PNP equations in ion channels using RFR. In response to the inefficiency and instability of traditional manual parameter selection, a Gummel iteration based on RFR parameter selection Algorithm is proposed. The algorithm addresses two major challenges in training data generation: first, an artificial ion channel generation algorithm is designed to efficiently construct channels with diverse geometric configurations, enriching the training dataset; second, two point cloud-based feature extraction algorithms are developed to effectively capture critical structural characteristics of the channels, providing physically meaningful input representations for the model. We conducted numerical experiments on three real ion channels (GA, VDAC and KcSA) with different geometric shapes. The results indicate that even when facing channel shapes different from the training set, The RFR based method for adaptive parameter selection not only significantly improves the convergence stability of Gummel iteration, but also outperforms traditional empirical parameter selection methods in terms of overall computational efficiency. Although the relaxation parameters manually selected are simple to choose based on experience, the problem is that when the same relaxation parameters are selected for different physical quantities (such as voltage and concentration), the number of iterations may vary greatly. At this time, selection based solely on experience becomes very unstable. By comparison, for the relaxation parameters selected by machine learning used in this paper, whether from the perspective of voltage or concentration, the number of iterations will be more stable compared with those selected based on experience. In addition, the proposed method can be extended to parameter selection in other multiphysics-coupled iterative numerical processes, while the training set generation strategy and feature extraction techniques also offer reliable data construction solutions for scientific computing problems involving irregular interfaces or complex geometries. \\

\noindent
{\bf Acknowledgement } Y. Yang was partially supported by the Science and Technology Project of Guangxi (Guike AD25069086), was supported by the Guangxi Natural Science Foundation, China (Grant No. 2025GXNSFAA069683) and National Natural Science Foundation of China (Grant No. 12561072). S. Shu was supported by National Natural Science Foundation of China (Grant No. 12371373) and Science Challenge Project of China (Grant No. TZ2024009).

\bibliography{artpnp}
\bibliographystyle{unsrt}
\end{document}